\documentclass[bibother]{asl}  
\usepackage{amssymb}
\usepackage{graphicx}
\usepackage{framed}
\usepackage{tikz,subfigure}
\usepackage{amssymb}
\usepackage{graphicx}
\usepackage{graphics} 
\usepackage[all,arc,poly]{xy}
\usepackage{color}
\usepackage{relsize}

\renewcommand{\phi}{\varphi}

\newcommand\<{(}
\renewcommand\>{)}

\newcommand{\nc}{\newcommand}
\nc{\look}{\marginpar{$\bullet$}}
\nc{\Section}{\section}
\nc{\SubSection}{\subsection}

\newtheorem{theo}{Theorem}[section]   
\newtheorem{ddef}[theo]{Definition}
\newtheorem{llem}[theo]{Lemma} 
\newtheorem{oobs}[theo]{Observation} 
\newtheorem{rrem}[theo]{Remark} 
\newtheorem{prop}[theo]{Proposition} 
\newtheorem{ccor}[theo]{Corollary}  
\newtheorem{qquest}[theo]{Question} 
\newtheorem{fact}[theo]{Fact} 
\newtheorem{pprov}[theo]{Proviso}
\newtheorem{eexam}[theo]{Example} 

\nc{\bT}{\begin{theo}} 
\nc{\eT}{\end{theo}}
\nc{\bD}{\begin{ddef} \rm }
\nc{\eD}{\end{ddef}}
\nc{\bC}{\begin{ccor}}
\nc{\eC}{\end{ccor}}
\nc{\bCl}{\begin{claim}}
\nc{\eCl}{\end{claim}}
\nc{\bQ}{\begin{qquest}}
\nc{\eQ}{\end{qquest}}
\nc{\bL}{\begin{llem}}
\nc{\eL}{\end{llem}}
\nc{\bP}{\begin{prop}}
\nc{\eP}{\end{prop}}
\nc{\bR}{\begin{rrem}}
\nc{\eR}{\end{rrem}}
\nc{\bO}{\begin{oobs}}
\nc{\eO}{\end{oobs}}
\nc{\bF}{\begin{fact}}
\nc{\eF}{\end{fact}}
\nc{\bProv}{\begin{pprov}}
\nc{\eProv}{\end{pprov}}
\nc{\bE}{\begin{eexam} \rm }
\nc{\eE}{\end{eexam}}

\renewcommand{\geq}{\geqslant}
\renewcommand{\leq}{\leqslant}

\renewcommand{\subset}{\subseteq}

\renewcommand{\supset}{\supseteq}
\newcommand{\kbox}{\ensuremath{\square}}
\newcommand{\ibox}{\boxplus}
\def\disji{\rotatebox[origin=c]{-90}{$\!{\geqslant}$}}
\newcommand{\lori}{\,\disji\,}
\newcommand{\simn}{{\sim^n}}

\def\Disji{\rotatebox[origin=c]{-90}{$\!{\mathlarger{\mathlarger{\,\mathlarger{\geqslant}}}}\!$}}
\newcommand{\Lori}{\,\Disji\,} 
\newcommand{\ee}{\varepsilon}
\newcommand{\E}{\texttt{E}}
\newcommand{\eee}{\ensuremath{\epsilon}}
\newcommand{\PP}{\texttt{P}}

\nc{\simg}{\sim_{\mathsf{g}}}
\nc{\ssimg}{\approx_{\mathsf{g}}}
\nc{\FO}{\mathsf{FO}}
\nc{\MSO}{\mathsf{MSO}}
\nc{\ML}{\mathsf{ML}}
\nc{\IML}{\mathsf{IML}}
\nc{\SO}{\mathsf{SO}}
\nc{\GF}{\mathsf{GF}}
\nc{\gd}{\mathsf{gd}}
\nc{\hg}{\mathsf{hg}}
\nc{\vg}{\mathsf{vg}}
\nc{\free}{\mathrm{free}}
\nc{\dom}{\mathrm{dom}}
\nc{\nn}{\mathsf{non}}
\nc{\NE}{\mathsf{NE}}
\nc{\NF}{\mathsf{NF}}
\nc{\KK}{\ensuremath{\mathcal{K}}}

\nc{\C}{\mathrm{cl}}

\nc{\full}{\mathsf{full}}
\nc{\lf}{\mathsf{lf}}
\nc{\dc}{\mathsf{dc}}
\nc{\rel}{\mathsf{rel}}

\newenvironment{romanenumerate}%
{\begin{list}{(\roman{enumi})}{\usecounter{enumi}
\setlength{\labelwidth}{2cm}
\setlength{\itemindent}{0pt}
\setlength{\itemsep}{0.5\itemsep}
\setlength{\topsep}{\itemsep}
\setlength{\parsep}{0pt}
}}{\end{list}}
\nc{\bre}{\begin{romanenumerate}}
\nc{\ere}{\end{romanenumerate}}
\newenvironment{alphaenumerate}%
{\begin{list}{(\alph{enumii})}{\usecounter{enumii}
\setlength{\labelwidth}{2cm}
\setlength{\itemindent}{0pt}
\setlength{\itemsep}{0.5\itemsep}
\setlength{\topsep}{\itemsep}
\setlength{\parsep}{0pt}
}}{\end{list}}
\nc{\bae}{\begin{alphaenumerate}}
\nc{\eae}{\end{alphaenumerate}}
\newenvironment{numenumerate}%
{\begin{list}{(\arabic{enumiii})}{\usecounter{enumiii}
\setlength{\labelwidth}{2cm}
\setlength{\itemindent}{0pt}
\setlength{\itemsep}{0.5\itemsep}
\setlength{\topsep}{\itemsep}
\setlength{\parsep}{0pt}
}}{\end{list}}
\nc{\bne}{\begin{numenumerate}}
\nc{\ene}{\end{numenumerate}}

\nc{\ins}[1]{\bigskip\noindent
\framebox{\begin{minipage}{.98\textwidth} \sloppy \noindent \em #1 \end{minipage}}\bigskip}

\nc{\str}[1]{{\mathfrak{#1}}}
\nc{\brck}[1]{[\![ #1 ]\!]}
\nc{\restr}{\!\restriction\!}

\nc{\HH}{\mathbb{H}}
\nc{\VV}{\mathbb{V}}
\nc{\abar}{\mathbf{a}}
\nc{\bbar}{\mathbf{b}}
\nc{\cbar}{\mathbf{c}}
\nc{\xbar}{\mathbf{x}}
\nc{\ybar}{\mathbf{y}}
\nc{\zbar}{\mathbf{z}}
\nc{\ubar}{\mathbf{u}}
\nc{\sbar}{\mathbf{s}}
\nc{\tbar}{\mathbf{t}}
\nc{\vbar}{\mathbf{v}}
\nc{\wbar}{\mathbf{w}}

\nc{\Xbar}{\mathbf{X}}
\nc{\Ybar}{\mathbf{Y}}
\nc{\Zbar}{\mathbf{Z}}
\nc{\Pbar}{\mathbf{P}}
\nc{\nubar}{\mbox{\boldmath $\nu$}}

\nc{\barr}{\begin{array}}
\nc{\earr}{\end{array}}
\nc{\btab}{\begin{tabular}}
\nc{\etab}{\end{tabular}}

\nc{\nothing}{\rule{0em}{1ex}}
\nc{\highnothing}{\rule{0em}{3ex}}
\nc{\hnt}{\highnothing}
\nc{\nt}{\nothing}
\nc{\nnt}{\rule{.1pt}{0pt}}

\nc{\ssc}{\scriptscriptstyle}
\nc{\N}{{\mathbb N}}
\nc{\Z}{{\mathbb Z}}
\nc{\M}{{\mathbb M}}
\nc{\F}{{\mathbb F}}
\nc{\W}{{\mathbb W}}
\newcommand{\PlayerI}{\mbox{\bf I}}
\newcommand{\PlayerII}{\mbox{\bf II}}
\renewcommand{\P}{\ensuremath{\mathcal{P}}}
\newcommand{\A}{\ensuremath{\mathcal{A}}}
\newcommand{\CC}{\ensuremath{\mathcal{C}}}

\newcommand{\MM}{\ensuremath{\mathfrak{M}}}
\newcommand{\stw}{\mathrm{ST}_{\mathtt{w}}}

\newcommand{\sts}{\mathrm{ST}_{\mathtt{s}}}
\newcommand{\stt}{\mathrm{ST}_{\mathtt{t}}}
\newcommand{\e}{\texttt{e}}
\newcommand{\s}{\mathtt{s}}

\renewcommand{\t}{\mathtt{t}}
\newcommand{\w}{\mathtt{w}}
\renewcommand{\v}{\mathtt{v}}

\newcommand{\inqbm}{\textsc{InqML}}
\newcommand{\inqml}{\textsc{InqML}}

\nc{\prf}{\begin{proof}}
\nc{\eprf}{\end{proof}}

\usepackage{wrapfig}

\tikzstyle{index on}=[inner sep=2pt, white, circle, fill=black]
\tikzstyle{index off}=[inner sep=2pt, black, circle, draw]
\tikzstyle{index gray}=[inner sep=2pt, black, circle, fill=lightgray]
\tikzstyle{opaque}=[fill=gray,fill opacity=.1]

\title{Inquisitive Bisimulation} 

\author{Ivano Ciardelli and Martin Otto}

\thanks{%
The present paper extends and subsumes the proceedings paper \cite{ICMO17}. It is a revised version of the first part of a larger manuscript (\cite{CiardelliOttoarXiv18}), accessible via arXiv, the second part of which we intend to develop into a separate paper. We are grateful to the reviewers for their precious comments and suggestions, which pushed us to clarify some important points.
Martin Otto's research was
    partially funded by DFG grant OT~147/6-1: \emph{Constructions and
      Analysis in Hypergraphs of Controlled Acyclicity}; it was also
    greatly supported through his participation 
in a programme on \emph{Logical Structure in Computation} at the
Simons Institute in Berkeley in 2016, which is gratefully
acknowledged. 
}

\begin{document}
\maketitle
 
\begin{abstract}
\noindent
Inquisitive modal logic, \inqml, is a generalisation of standard Kripke-style 
modal logic. In its epistemic incarnation, it extends standard epistemic logic 
to capture not just the information that agents have, but also the questions 
that they are interested in. Technically, \inqml\ fits within the family of 
logics based on team semantics. From a model-theoretic perspective, it 
takes us a step in the direction of monadic second-order logic, as inquisitive 
modal operators involve quantification over sets of worlds. 
We introduce and investigate the natural notion of bisimulation
equivalence  in the setting of \inqml. 
We compare the expressiveness of \inqml\ and first-order logic 
in the context of relational structures with two sorts, one 
for worlds and one for information states, and characterise  
inquisitive modal logic 
as the bisimulation invariant fragment of first-order logic over various 
natural classes of two-sorted structures. 
\end{abstract}

\section{Introduction}

The recently developed framework of \emph{inquisitive logic}
\cite{CiardelliRoelofsen,Ciardelli:15inqd,IvanoDiss,Ciardelli:18qait}
can be seen as a generalisation of classical logic which encompasses
not only statements, but also questions. One reason why this
generalisation is interesting is that it provides a novel perspective
on the logical notion of \emph{dependency}, which plays an important
r\^{o}le in applications (e.g., in database theory) and which has
recently received attention in the field of \emph{dependence logic}
\cite{Vaananen:07}.  Indeed, dependency is nothing but a facet of the
fundamental logical relation of entailment, once this is extended so
as to apply not only to statements, but also to
questions~\cite{Ciardelli:16dependency}. This connection explains the
deep similarities existing between systems of inquisitive logic and
systems of dependence logic (see~\cite{Yang:14,Ciardelli:16dependency,IvanoDiss,YangVaananen}).
A different r\^{o}le for questions in a logical system comes from the
setting of modal logic: once the notion of a modal operator is
suitably generalised, questions can be embedded under modal operators
to produce new statements that have no ``standard'' counterpart. This
approach was first developed in~\cite{CiardelliRoelofsen:15idel} in
the setting of epistemic logic. The resulting \emph{inquisitive
epistemic logic} 
models not only the information that agents
have, but also the issues that they are interested in, i.e., the
information that they would like to obtain. Modal formulae in
inquisitive epistemic logic
can express not only that an agent knows that~$p$ (by the formula $\kbox p$) but also
that she knows \emph{whether}~$p$ ($\kbox{?p}$) or that she
\emph{wonders whether}~$p$ (by the formula $\ibox{?p}$)---a statement that cannot be
expressed without the use of embedded questions. As shown 
in~\cite{CiardelliRoelofsen:15idel}, several key notions of epistemic
logic generalise smoothly to questions: besides common knowledge we
now have \emph{common issues}, the issues publicly entertained by the
group; and besides publicly announcing a statement, agents can now
also publicly ask a question, which typically results in new common
issues.  Thus, inquisitive epistemic logic may 
be seen as one step in extending modal logic
from a framework to reason about information and information change,
to a richer framework which also represents a higher stratum of
cognitive phenomena, in particular issues that may be raised in a
communication scenario.

Of course, like standard modal logic, inquisitive modal logic provides
a general framework that admits various interpretations, each
suggesting corresponding constraints on models.  For instance, an
interpretation of \inqml\ as a logic of action is suggested 
in~\cite{IvanoDiss}. In that interpretation, a modal formula $\kbox{?p}$
expresses that whether a certain fact $p$ will come about is
pre-determined independently of the agent's choices, while $\ibox{?p}$
expresses that whether $p$ will come about is fully determined by the
agent's choices.

From the perspective of mathematical logic, inquisitive modal logic is a natural
generalisation of standard modal logic. In standard modal logic, the accessibility
relation of a Kripke model associates 
with each possible world $w\in W$ a set $\sigma(w) \subset W$ of
possible worlds, namely, the worlds accessible from $w$; 
any formula $\phi$ of modal logic is semantically
associated with a set $|\phi|_\M \subset W$ of worlds, namely, the set of worlds where it is true; modalities then express relationships between these sets: 
for instance, $\kbox\phi$ expresses the fact that $\sigma(w)\subseteq|\phi|_\M$.  
In the inquisitive setting, the situation is similar: we still have a set $\Sigma(w)$ associated with each possible world $w$, and a set $[\phi]_\M$ associated with a formula $\phi$. Now, however, both $\Sigma(w)$ and $[\phi]_\M$ are no longer sets of worlds, but sets \emph{of sets} of
worlds. Inquisitive modalities still express
relationships between 
these two objects: $\kbox\phi$ expresses the fact that
$\bigcup\Sigma(w)\in[\phi]_\M$, while $\ibox\phi$ 
expresses the fact that $\Sigma(w)\subseteq[\phi]_\M$.

In this manner, inquisitive logic leads to a new framework for modal
logic that can be viewed as a generalisation of the standard framework. 
This raises the question of whether and how the classical
notions and results of modal logic carry over to this more general
setting. In this paper we address this question for the fundamental notion 
of \emph{bisimulation} and for two classical results revolving around
this notion, namely, the Ehrenfeucht-Fra\"iss\'e theorem for modal 
logic, and van Benthem style characterisation
theorems~\cite{GorankoOtto,Benthem83,Rosen,OttoNote}. 
A central~topic of this paper is the r\^ole of \emph{bisimulation invariance} 
as a unifying semantic feature that distinguishes modal logics
from classical predicate logics. As in many other areas, from 
temporal and process logics to knowledge representation in AI 
and database applications, so also in the inquisitive setting 
we find that the appropriate notion of  bisimulation invariance allows for 
precise model-theoretic characterisations of the expressive power of modal logic 
in relation to first-order logic.

Our first result is that the right notion of
inquisitive bisimulation equivalence $\sim$, with finitary approximation levels 
$\simn$, supports a counterpart of 
the classical  Ehrenfeucht--Fra\"\i ss\'e characterisations 
for first-order logic or for basic modal logic. 
This result establishes an exact correspondence between the expressive
power of \inqml\ and the \emph{finite} approximation
levels of inquisitive bisimulation equivalence:
if two points are behaviourally different in a way that can be detected within a finite 
number of steps, then the difference between them is witnessed by an
\inqml\ formula, and vice versa.
The result is non-trivial in our setting 
because of some subtle issues stemming from the interleaving of first- and
second-order features in inquisitive modal logic. 

\medskip
\bT[inquisitive Ehrenfeucht--Fra\"\i ss\'e theorem]\label{EFthm}~\\
Over finite vocabularies, the finite levels $\simn$ of 
inquisitive bisimulation equivalence correspond to 
the levels of $\inqbm$-equivalence up to modal nesting depth~$n$. 
\label{mainEF}
\eT

In order to compare \inqbm\ with classical first-order logic, we define a class
of two-sorted relational structures, and show how such structures encode 
models for \inqbm.  With respect to such relational structures
we find not only a ``standard translation'' of \inqbm\   
into two-sorted first-order logic, but also a van Benthem style
characterisation of $\inqbm$ as the bisimulation-invariant fragment of 
(two-sorted) first-order logic over several natural classes of
models. These results are technically
interesting, and they are not available on the basis of classical techniques,
because the relevant classes of two-sorted models are 
non-elementary (in fact, first-order logic is not compact over these classes, as we show). 
Our techniques yield characterisation theorems both in the setting of arbitrary 
 inquisitive models, and in restriction to just finite ones---i.e.\
 both in the sense of classical model theory and in the sense of
 finite model theory.

\bT
\label{main1}
Inquisitive modal logic can be characterised as the 
$\sim$-invariant fragment of first-order logic $\FO$ over natural 
classes of (finite or arbitrary) relational inquisitive models. 
\eT

Beside the conceptual development and the core results themselves, we think that 
also the methodological aspects of the present investigations have
some intrinsic value. Just as inquisitive logic models cognitive phenomena  
at a level strictly above that of standard modal logic, so the
model-theoretic analysis moves up from the level of ordinary first-order logic 
to a level strictly between first- and second-order logic. This level
is realised by first-order logic in a two-sorted framework that
incorporates second-order objects in the second sort in a controlled fashion. 
This leads us to generalise a number of notions and
techniques developed in the model-theoretic analysis of modal logics 
over non-elementary classes of frames 
(cf.~\cite{GorankoOtto,OttoNote,DawarOttoAPAL09,OttoAPAL04}, among others).
In the present paper we technically focus on the general case of 
inquisitive modal models. This also sets the stage for the model-theoretic 
treatment of inquisitive epistemic models. That case, which is of 
particular interest from the point of view of logical modelling, also requires 
some further extensions of the technical apparatus presented here. 
We aim to present corresponding results from~\cite{CiardelliOttoarXiv18}
in a sequel to the present paper.

\section{Inquisitive modal logic}
\label{sec:inquisitive modal logic}

In this section we provide an essential introduction to inquisitive
modal logic, \inqbm\ \cite{IvanoDiss}.
For further details and proofs, we refer to \S7 of~\cite{IvanoDiss}.

\subsection{Foundations of inquisitive semantics}

Usually, the semantics 
of a logic specifies truth-conditions for the formulae of the logic. 
In modal logics these truth-conditions are relative to possible worlds
in a Kripke model. 
However, this approach is limited in an important way: while suitable for statements, 
it is inadequate for questions. To overcome this limitation, inquisitive logic
interprets formulae not relative to states of affairs (possible worlds), but
relative to states of information. Following a tradition that goes
back to the work of Hintikka~\cite{Hintikka:62},  information states
are modelled extensionally as sets of worlds, namely, the set of those worlds
which are compatible with the given~information.\footnote{An analogous 
step from single worlds to sets of worlds (or, depending on the 
setting, from assignments to sets of assignments) is taken in recent 
 work on independence-friendly logic~\cite{Hodges:97,Hodges:97b} 
and dependence logic~\cite{Vaananen:07,Vaananen:08,AbramskyVaananen:09,Galliani:12,Yang:14,YangVaananen}, 
where sets of worlds are referred to as \emph{teams}. 
Although they originated independently and for different purposes, 
inquisitive logic and dependence logic are tightly related. For
detailed discussion of this connection, 
see~\cite{IvanoDiss,Ciardelli:16dependency}.\label{fn:team}}

\bD[information states]~\\  
An \emph{information state} over a set of worlds $W$ is a subset
$s\subseteq W$.
\eD

The empty set represents a state of inconsistent information, 
which is not compatible with any world. 
We refer to it as the \emph{inconsistent state}.
 
Rather than specifying when a sentence is \emph{true} at a world $w$, 
inquisitive semantics specifies when a sentence is \emph{supported} by an
information state $s$: intuitively, for a statement $\alpha$ this means that 
the information available in
$s$ implies that $\alpha$ is true; for a question~$\mu$, it means that
the information available in $s$ settles~$\mu$.
If $t$ and $s$ are information states and $t\subseteq s$, this means that
$t$ holds at least as much information as $s$: we say that $t$ is an 
\emph{extension} of~$s$. If $t$ is an
extension of $s$, everything that is supported at $s$ will also
be supported at $t$. This is a key feature of inquisitive
semantics, and it leads naturally to the notion of an \emph{inquisitive state}.

\bD[inquisitive states]~\\
An \emph{inquisitive state} over a set of possible worlds 
$W$ is a 
non-empty set of information states $\Pi\subseteq\wp(W)$ that is
\emph{downward closed} in the sense that 
$s\in\Pi$ implies  $t \in\Pi$ for all $t\subseteq s$ (downward 
closure).
\eD

The downward closure condition requires that $\Pi$ be closed under extensions of  information states. 
As described in the next section, an inquisitive state can be seen as a combined representation of information and issues. For more discussion on the significance of this structure, see~\cite{Ciardelli:13compass,CiardelliRoelofsen:15idel,Ciardelli:18book}.

\subsection{Inquisitive modal models}

A Kripke frame can be thought of as a set $W$ of worlds 
together with a map $\sigma$ that equips each world with a set of 
worlds $\sigma(w)$, i.e., an information state: 
the set of worlds that are \emph{accessible} from~$w$. 

Similarly, an inquisitive modal frame consists of a set $W$ of worlds
together with an \emph{inquisitive assignment}, a map $\Sigma:W\to\wp\wp(W)$
that assigns to each world a corresponding inquisitive state
$\Sigma(w)$, i.e.\ a downward closed set of information 
states. 
A \emph{model} is a frame enriched by a 
propositional assignment.

\bD[inquisitive modal models]~\\
An inquisitive modal frame is a pair $\F=(W,\Sigma)$, where $W$ is a
set, whose elements are referred to as \emph{worlds}, and
$\Sigma\colon W\to\wp\wp(W)$ assigns 
an inquisitive state $\Sigma(w)$ to each world $w\in W$.

An inquisitive modal model for a set \P\ of propositional atoms is a pair
$\M=(\F,V)$ where $\F$ is an inquisitive modal frame, and 
$V \colon \P \rightarrow \wp(W)$
is a propositional assignment.

A {world-(or state-)pointed} inquisitive modal model is a pair
consisting of a model $\M$ and a distinguished world (or state) in $\M$.
\eD

With an inquisitive modal model $\M$ we can always associate a
standard Kripke model $\str{K}(\M)$ having the same set of worlds and  
modal accessibility map $\sigma:W\to\wp(W)$ induced by the inquisitive
map $\Sigma$ according to 
\[
\barr{rcl}
\sigma \,\colon\, W &\longrightarrow&\wp(W)
\\
w &\longmapsto& \sigma(w):=\bigcup\Sigma(w).
\earr
\]

A natural interpretation for inquisitive modal models is the epistemic
one, developed in~\cite{CiardelliRoelofsen:15idel,Ciardelli:14aiml}. 
In that interpretation, the map $\Sigma$ is taken to describe not only 
an agent's \emph{knowledge}, as in standard epistemic logic, but also 
an agent's \emph{issues}.%
\footnote{In inquisitive semantics, the term \emph{issue} is used to refer 
to the content of a question. For instance, the issues that a detective 
entertains might be those expressed by the questions 
\emph{who committed the murder}, \emph{whether they had an accomplice}, 
and \emph{what the motive is}.} 
The agent's knowledge state at $w$, $\sigma(w)=\bigcup\Sigma(w)$, 
consists of all the worlds that are compatible with what the agent knows. 
The agent's inquisitive state at $w$, $\Sigma(w)$, consists of all those 
information states where the agent's issues are settled. This interpretation 
is particularly interesting in the multi-modal setting, where a model comes 
with multiple state maps $\Sigma_a$, one for each agent $a$ in a set~$\A$. 
Moreover, this specific interpretation suggests some constraints on the maps 
$\Sigma_a$, analogous to the usual $S5$ constraints on Kripke models.

\bD[inquisitive epistemic models]
\label{epinqdef}\\
An \emph{inquisitive epistemic frame} for a set $\A$ of agents is a pair
$\F=( W,(\Sigma_a)_{a\in\A})$, where each map
$\Sigma_a:W\to\wp\wp(W)$ 
assigns to each world
$w$ an inquisitive state $\Sigma_a(w)$ in accordance with the 
following constraints, where $\sigma_a(w)=\bigcup\Sigma_a(w)$:\!\!\nt
\bre
\item[--] $w\in\sigma_a(w)$ (factivity);
\item[--] $v\in\sigma_a(w) \;\Rightarrow\; \Sigma_a(v)=\Sigma_a(w)$
(full introspection).
\ere
\eD

It is easy to verify that the Kripke frame associated with an inquisitive
epistemic frame is an $S5$ frame, i.e., the accessibility
maps $\sigma_a$ correspond to accessibility relations 
$R_a \!:=\! \{ ( v,w) \colon \!v\! \in\! \sigma_a(w) \}$
that are equivalence relations on~$W$.

\bE\label{ex:1}
Consider a model with four worlds, $w_{pq},w_{p\overline
  q},w_{\overline p q},w_{\overline p\overline q}$, where the
subscript indicate the propositional valuation at each world. The
inquisitive state map $\Sigma$ is as follows, where $S^\downarrow$
indicates the closure of the set $S \subset \wp(W)$ 
under subsets.
\[
\barr{@{}l@{}}
\Sigma(w_{pq})=\Sigma(w_{p\overline{q}})=\{\{w_{pq}\},\{w_{q\overline{q}}\}\}^\downarrow
\\
\hnt
\Sigma(w_{\overline pq})=\Sigma(w_{\overline p\overline
  q})=\{\{w_{\overline pq},w_{\overline p\overline q}\}\}^\downarrow
\earr 
\]
\begin{figure}
  \begin{center}
   \begin{tikzpicture}[>=latex,scale=0.6]
 \draw[opaque,rounded corners] (-1.7,1.7) rectangle (-.3, .3);
 \draw[opaque,rounded corners] (1.7,1.7) rectangle (.3, .3);
 \draw[rounded corners,dashed] (-1.85,1.85) rectangle (1.85, .15);
 \draw[opaque,rounded corners,yshift=-2cm] (-1.7,1.7) rectangle (1.7, .3);
 \draw[rounded corners,dashed,yshift=-2cm] (-1.85,1.85) rectangle (1.85, .15);
 \draw (1.9,1.9) rectangle (1.9,1.9); 
 \draw (-1.9,-1.9) rectangle (-1.9,-1.9); 

 \draw (-1,1) node[index gray] (yy) {$w_{pq}$};
 \draw (1,1) node[index gray] (yn) {$w_{p\overline q}$};
 \draw (-1,-1) node[index gray] (ny) {$w_{\overline p q}$};
 \draw (1,-1) node[index gray] (nn) {$w_{\overline p\overline q}$};

	\end{tikzpicture}
	\caption{A single-agent inquisitive epistemic model}
	\label{fig:1}
\end{center}
\end{figure}
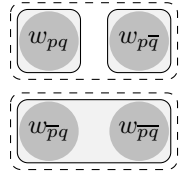 
This model is depicted in Figure~\ref{fig:1}. At a world $w$, the
epistemic state $\sigma(w)$ of the agent consists of those worlds
included in the same dashed area as $w$; the solid blocks inside this
area are the maximal elements of the inquisitive state $\Sigma(w)$---i.e., the maximal states in which the issue is resolved.

At worlds $w_{pq}$ and $w_{p\overline q}$, the agent's knowledge state is $\{w_{pq},w_{p\overline q}\}$: that is, the agent knows that $p$ is true, but not whether $q$ is true. Moreover, in order to settle the agent's issues it is necessary and sufficient to reach an extension of the current state which settles whether $q$. In short, then, these are worlds where the agent knows that $p$ and wonders whether $q$.

At worlds $w_{\overline pq}$ and $w_{\overline p\overline q}$, the agent's knowledge state is $\{w_{\overline pq},w_{\overline p\overline q}\}$: that is, the agent knows that $\neg p$, but not whether $q$. However, at these worlds no further information is needed to resolve the agent's issues. Thus, these are worlds where the agent knows that $\neg p$ and does not have any remaining issues.
\eE

\subsection{Inquisitive modal logic}

The syntax  of inquisitive modal logic $\inqbm$ is given by:
$$\phi::= p\;|\,\bot\,|\,(\phi\land\phi)\,|\,(\phi\to\phi)\,|\,(\phi\lori\phi)\,|\,\kbox\phi\,|\,\ibox\!\phi$$

The syntax of inquisitive epistemic logic is defined analogously,
except that modalities are indexed by agents; that is, for every agent
$a\in\A$ we have two corresponding modalities $\kbox_a$ and $\ibox_a$,
which are interpreted based on the state map $\Sigma_a$ associated
with the agent.\footnote{In~\cite{CiardelliRoelofsen:15idel,Ciardelli:14aiml} the modalities
$\kbox_a$ and $\ibox_a$ are denoted $K_a$ and $E_a$, and
read as ``know'' and ``entertain'' respectively.}
We treat negation and disjunction as defined connectives (syntactic
shorthands) according to
\[
\neg\phi:=\phi\to\bot
\qquad \mbox{ and } \qquad
\phi\lor\psi:=\neg(\neg\phi\land\neg\psi)
\]
so that 
the above syntax emulates standard propositional formulae
in terms of atoms and connectives $\land$ and $\to$ together with the defined
$\neg$ and $\vee$. The semantics of these 
will be essentially the same as in standard propositional logic. 

In addition to standard connectives, our language contains a new
connective, $\lori$, called \emph{inquisitive disjunction}. We may read
formulae built up by means of this connective as propositional
questions. For instance, we read the formula $p\lori\neg p$ as the question
\emph{whether or not $p$}, and we abbreviate this formula as $?p$.  
Our language also contains two modalities, which are allowed to
embed both statements and questions. 
As we shall see, both these modalities coincide with a standard Kripke
box modality when applied to statements, but crucially differ when applied to questions. 
In particular, under an epistemic interpretation $\kbox{?p}$ expresses 
the fact that the agent knows whether $p$, while $\ibox{?p}$ 
intuitively says, roughly, 
that the agent is interested in the issue whether~$p$. 

As mentioned above, the semantics of \inqbm\ is given in terms of
support in an information state, rather than truth at a
possible world.%
\footnote{This means that \inqbm\ fits within the quickly growing
  family of logics based on a \emph{team semantics}. See footnote~\ref{fn:team} on
  page~\pageref{fn:team} and the references therein.}

\bD[semantics of $\inqbm$]~\\
\label{supportsemdefn}%
Let $\M=(W,\Sigma,V)$ be an inquisitive modal model, $s\subseteq W$:
\begin{itemize}
\item $\M,s\models p\iff s\subseteq V(p)$
\item $\M,s\models \bot\iff s=\emptyset$
\item $\M,s\models \phi\land\psi\iff \M,s\models\phi\text{ and }\M,s\models\psi$
\item 
$\M,s\models \phi\to\psi\iff \forall t\subseteq s:\M,t\models\phi\Rightarrow \M,t\models\psi$
\item $\M,s\models \phi\lori\psi\iff \M,s\models\phi\text{ or }\M,s\models\psi$
\item $\M,s\models \kbox\phi\iff\forall w\in s: \M,\sigma(w)\models\phi$
\item $\M,s\models \ibox\phi\iff\forall w\in s\;\forall t\in\Sigma(w): \M,t\models\phi$
\end{itemize}
\eD

If a state $s$ can be extended consistently to a state that supports a formula $\phi$, we say that $s$ is \emph{compatible} with $\phi$:
\[
s \mbox{ is compatible with } \phi
\quad \mbox{ iff } \quad 
\exists t\subseteq s \colon t\neq\emptyset \mbox{ and }
\M,t\models\phi,
\]

The derived clauses for the defined connectives $\neg$ and $\lor$ then read as follows:
\begin{itemize}
\item $\M,s\models\neg\phi\iff s$ is not compatible with $\phi$
\item $\M,s\models\phi\lor\psi\iff\forall t\subseteq s, t\neq\emptyset: t$ is compatible with $\phi$ or with $\psi$
\end{itemize}

As an illustration, consider the support conditions for the formula
\[
?p:=p\lori\neg p.
\] 

This formula is supported by a state $s$ in case
$p$ is true at all worlds in $s$ (i.e., if the information available
in $s$ implies that $p$ is true) or in case $p$ is false at all worlds
in $s$ (i.e., if the information available in $s$ implies that $p$ is
false).  Thus, $?p$ is supported by those information states
that settle whether or not $p$ is true.

\bP\label{prop:persistency}
The following properties hold generally in \inqbm:
\begin{itemize}
\item persistency: if $\M,s\models\phi$ and $t\subseteq s$, then $\M,t\models\phi$;
\item semantic ex-falso: $\M,\emptyset\models\phi$ for all $\phi\in\inqbm$.
\end{itemize}
\eP

The first principle says that support is preserved as information
increases, i.e., as we move from a state to an extension of it. The
second principle says that the empty set of worlds---the
inconsistent 
information state---vacuously supports every formula.  Together, these
principles imply that the support set $[\phi]_\M := \{ s \subseteq W
\colon \M,s \models \phi \}$ of a formula is downward closed and
non-empty, i.e., it is an inquisitive state.

Although the primary notion of our semantics is support at an information state, truth at a world is retrieved by defining it as support with respect to singleton states.

\bD[truth]
\label{truthdef}~\\ 
$\phi$ is \emph{true} at a world $w$ in a model $\M$, 
denoted $\M,w\models\phi$, in case $\M,\{w\}\models\phi$.
\eD

Spelling out Definition~\ref{supportsemdefn} in the special case of
singleton states, we see that standard connectives have the usual
truth-conditional behaviour. 
For modal formulae, we find the following truth-conditions:

\begin{itemize}
\item $\M,w\models \kbox\phi\iff \M,\sigma(w)\models\phi$
\item $\M,w\models \ibox\phi\iff \forall t\in\Sigma(w): \M,t\models\phi$
\end{itemize}

Notice that truth in \inqbm\ cannot be given a direct recursive
definition, as the truth conditions for modal formulae $\kbox\phi$ and
$\ibox\phi$ depend on the support conditions for $\phi$---not just on
its truth conditions.

For many formulae, support at a state just boils down to truth at each
world. We refer to these formulae as \emph{truth-conditional}.%
\footnote{In team semantic terminology 
(e.g., \cite{Vaananen:07,YangVaananen}),
truth-conditional formulae are called \emph{flat}.}

 \bD[truth-conditional formulae]\label{flatdef}~\\
 We say that a formula $\phi$ is \emph{truth-conditional} if for all
 models $\M$ and information states $s$: $\M,s\models\phi\iff \M,w\models\phi$ for all $w\in s$.
 \eD

Following~\cite{IvanoDiss}, we view truth-conditional 
formulae as \emph{statements}, and non-truth-conditional formulae as
\emph{questions}. 
The next proposition identifies a large class of formulae that are truth-conditional.

\bP\label{prop:truth-conditionality} \mbox{}
\label{flatprop}
Atomic formulae (including $\bot$) 
and all formulae of the form $\kbox\phi$ and $\ibox\phi$ are truth-conditional.
The class of truth-conditional formulae is closed under all
connectives except for $\lori$.
\eP

Using this fact, it is easy to see that all formulae of standard modal logic, i.e., formulae which do not contain $\lori$ or $\ibox$, receive exactly the same truth-conditions as in standard modal logic.

\bP 
If $\phi$ is a formula not containing $\lori$ or $\ibox$, then we have $\M,w\models\phi$ if and only if $\str{K}(\M),w\models\phi$ holds in standard Kripke semantics.
\eP

As long as questions are not around, the modality~$\ibox$ also
coincides with $\kbox$, and with the standard box modality. 
That is, if $\phi$ is truth-conditional, then 
\[
M,w\models\kbox\phi\iff M,w\models\ibox\phi\iff M,v\models\phi\text{
  for all }v\in\sigma(w).
\]

Thus, the two modalities coincide on statements. However, they come
apart when they are applied to questions. For an illustration,
consider the formulae $\kbox{?p}$ and $\ibox{?p}$ in the epistemic
setting: $\kbox{?p}$ is true iff the knowledge state of the agent,
$\sigma(w)$, settles the question $?p$; thus, $\kbox{?p}$ expresses
the fact that the agent knows whether $p$. By contrast, $\ibox{?p}$ is
true iff any information state $t\in\Sigma(w)$, i.e., any state that
settles the agent's issues, also settles $?p$; thus $\ibox{?p}$
expresses that finding out whether $p$ is part of the agent's goals.

\bE
\label{ex:2} 
Consider again the model of Example~\ref{ex:1}. 
The agent's knowledge state at world $w_{pq}$ is
$\sigma(w_{pq})=\{w_{pq},w_{p\overline q}\}$. Since
$\{w_{pq},w_{p\overline q}\}$ does not support~$?q$ we have
$\M,w\models\neg\kbox{?q}$. On the other hand, since the agent's
inquisitive state is $\Sigma(w_{pq})=\{\{w_{pq}\},\{w_{p\overline
q}\}\}^\downarrow$, and since each element in this state supports
$?q$, we do have $\M,w_{pq}\models\ibox{?q}$. This witnesses that, at
world $w_{pq}$, the agent does not know whether $q$ ($\neg\kbox?p$),
but she's interested in finding out ($\ibox?q$). By contrast, one can
check that at world $w_{\overline p q}$ we have $\M,w_{\overline p
q}\models \neg\kbox{?q}\land\neg\!\ibox\!{?q}$, witnessing that at
this world, the agent is neither informed about whether $q$, nor
interested in finding out.
\eE

\subsection{Defining properties of worlds and states}

Inquisitive modal formulae can be interpreted both relative to information states and (derivatively) relative to worlds. They can thus be seen both as expressing properties of state-pointed models, and as expressing properties of world-pointed models. We can identify these properties with the corresponding classes: 
 
\begin{itemize}
\item $\KK_\phi^{\w}=\{(M,w)\colon M,w\models\phi\}$
\item $\KK_\phi^{\s}=\{(M,s)\colon M,s\models\phi\}$
\end{itemize}

More generally, by a \emph{property} of world- or state-pointed models
we mean a class of such objects. We say that a property $\KK$ of
world-pointed models is definable in \inqml\ if $\KK=\KK_\phi^{\w}$ for
some formula $\phi$ of \inqml. Similarly, a property $\KK$ of
state-pointed models is definable in \inqml\ if $\KK=\KK_\phi^{\s}$ for
some $\phi$.

We can now formulate the main question that we will adress in this paper: 
which properties of world- or state-pointed models are definable in \inqml? 

For the case of state-pointed models, persistency and the semantic 
ex-falso condition (Proposition~\ref{prop:persistency}) impose an immediate 
constraint: in order for a property $\KK$ of state-pointed models to be 
definable in \inqml, $\KK$ must be an \emph{inquisitive property}, 
in the following sense.

\bD[inquisitive properties]\label{inquisitive state-properties}~\\
A property $\KK$ of state-pointed models is an \emph{inquisitive} 
property if the following two conditions hold:
\bre
\item if $(\M,s)\in\KK$ and $t\subseteq s$, then $(\M,t)\in\KK$;
\item $(\M,\emptyset)\in\KK$ for any model $\M$.
\ere
\eD

In the rest of the paper, when dealing with properties 
of state-pointed models, we can restrict our attention to inquisitive properties. 

What features must a world-property have in order to be \inqml\ definable? Similarly, what features must an inquisitive state-property have? The following two sections provide a precise answer to this question. 

\section{Inquisitive bisimulation}
\label{sec:bisimulations}

An inquisitive modal model can be seen as a structure with two sorts
of entities, worlds and information states, which interact with each
other. On one hand, an information state $s$ is completely
determined by the worlds that it contains; on the other hand, a world
$w$ is determined by the atoms it makes true and the information
states which lie in $\Sigma(w)$. Taking a more behavioural
perspective, we can look at an inquisitive modal model as a model
where two kinds of transitions are possible: from an information state
$s$ we can make a transition to a world $w\in s$, and from a world
$w$ we can make a transition to an information state $s\in\Sigma(w)$. 
This suggests a natural notion of bisimilarity, together with its natural
finite approximations of $n$-bisimilarity for $n \in \N$. As usual, these
notions can equivalently be defined either in terms of back-and-forth
systems or in terms of strategies in corresponding bisimulation games. 
We start from the latter due to its more immediate and intuitive appeal to the 
underlying dynamics of a ``probing'' of behavioural equivalence. 

The inquisitive bisimulation game is played by two players, \PlayerI\ and \PlayerII, 
who act as challenger and defender of a similarity
claim involving a pair of worlds $w$ and $w'$ or 
information states $s$ and $s'$ over two models
$\M=( W,\Sigma,V )$ and $\M'=( W',\Sigma',V')$.
We denote world-positions as $( w,w' )$ and 
state-positions as $( s,s' )$, where $w \in W, w' \in W'$
and $s \in \wp(W), s' \in \wp(W')$, respectively.
The game proceeds in rounds that alternate between world-positions and
state-positions. 
Playing from a world-position $( w,w' )$,
  \PlayerI\ chooses an information state in the inquisitive state associated to
 one of these worlds ($s \in \Sigma(w)$ or 
$s' \in \Sigma'(w')$)  
and \PlayerII\ must respond by choosing an information state on the
opposite side,
which results in a state-position $( s,s' )$.  
Playing from a state-position
$( s,s')$,
  \PlayerI\ chooses a world in either state ($w \in s$ or $w' \in s'$)
and \PlayerII\ must respond by choosing a world from the other state,
which results in a world-position $( w,w' )$.
A round of the game consists 
of four moves leading from a world-position to another.

In the bounded version of the game, the number of rounds is fixed
in advance. In the unbounded version, the game is allowed to go on
indefinitely. Either player loses when stuck for a move. The game ends
with a loss for \PlayerII\ in any world-position $( w,w')$ that shows a
discrepancy at the atomic level, i.e., such that $w$ and $w'$ disagree on the truth of some $p\in\P$. 
All other plays, and in particular infinite runs of the unbounded game, are won by $\PlayerII$.

\bD[bisimulation equivalence]\\
\label{bisimdef}%
Two world-pointed models $\M,w$ and $\M',w'$ are \emph{$n$-bisimilar}, 
$\M,w\,\simn\, \M',w'$, if \PlayerII\ has a winning strategy in the $n$-round 
game starting from $( w,w')$.
$\M,w$ and $\M',w'$ are \emph{bisimilar}, denoted 
$\M,w\sim \M',w'$, if \PlayerII\ has a winning strategy in the unbounded
 game starting from $( w,w')$. 
Two state-pointed models $\M,s$ and $\M',s'$ are ($n$-)bisimilar, denoted $\M,s\!\sim\!\M',s'$
(or $\M,s\,\simn\,\M',s'$), if every world in $s$
is ($n$-)bisimilar to some world in $s'$ and vice versa.%
\footnote{This definition of bisimilarity between  states is
  reminiscent of the corresponding definition given
  in~\cite{Kontinen:15} for modal team logic. Like inquisitive modal
  logic, modal team logic interprets formulae with respect to sets of
  possible worlds, and thus can be seen as expressing properties of
  state-pointed models. However, there are some major differences with
  the present setting. Most importantly, the structures for modal team
  logic are standard Kripke models. By contrast, \inqml\ is
  interpreted on models having a richer structure; information states
  enter the picture not just as evaluation points, but also in
  determining the structure of the model itself, since each world is
  associated with a set $\Sigma(w)$ of ``successors'' which are not
  worlds, but information states. 
This difference is reflected in the respective notions of bisimulation. In modal team logic, bisimilarity between worlds is the standard notion, and bisimilarity between states is a simple derivative of it: two states are bisimilar if each world in the one is bisimilar to some world in the other. By contrast, in our setting, world-bisimilarity and state-bisimilarity are inextricably intertwined. 
It is helpful to view this in terms of the bisimulation game. The game
for modal team logic starts with a pair of information states; player
\PlayerI\ selects a world from either state, and player \PlayerII\ responds with
a world in the other; after that, the standard bisimulation game for
modal logic is played. Thus, information states play a very limited
role: they only matter for the initial move, and moreover, there is no
move where players have to pick an information state. By contrast, in
the case of inquisitive modal logic, the game alternates indefinitely
between world positions and state positions, and moves in which
players pick information states are a crucial part of the game.}
\eD

These notions generalise naturally to the multi-modal setting 
with inquisitive assignments $(\Sigma_a)_{a \in \A}$ 
for a set $\A$ of agents; at a world-position, player \PlayerI\ also gets the choice of which agent to~probe.

Now let us turn to the static perspective on inquisitive bisimulations. One natural way to define a bisimulation between two models $\M$ and $\M'$ is as a relation which pairs up both the worlds and the states of these two models in such a way as to guarantee a winning strategy in the unbounded bisimulation game. This leads to the following definition.

\bD[bisimulation relations]~\\
Let $M=(W,\Sigma,V)$ and $M'=(W',\Sigma',V')$ be two inquisitive modal
models. A non-empty relation $Z\subseteq W\times W'\cup
\wp(W)\times\wp(W')$ is called a bisimulation in case the following 
constraints are satisfied:
\begin{itemize}
\item atom equivalence: if $wZw'$ then for all $p\in\P$, $w\in V(p)\iff w'\in V'(p)$
\item state-to-world back\&forth: if $sZs'$ then 
\begin{itemize}
\item for all $w\in s$ there is some $w'\in s'$ s.t.\ $wZw'$
\item for all $w'\in s'$ there is some $w\in s$ s.t.\ $wZw'$
\end{itemize}
\item world-to-state back\&forth: if $wZw'$ then 
\begin{itemize}
\item for all $s\in \Sigma(w)$ there is some $s'\in\Sigma'(w')$ s.t.\ $sZs'$
\item for all $s'\in \Sigma'(w')$ there is some $s\in\Sigma(w)$ s.t.\ $sZs'$
\end{itemize}
\end{itemize}
\eD

It is then routine to check that bisimilarity can be characterised in
terms of the existence of a bisimulation relation.

\bP 
Let $\M,x$ and $\M',x'$ be two world- or state-pointed models. 
$\M,x\sim\M',x'\iff$ there exists a bisimulation $Z$ such that $xZx'$. 
\eP

Alternatively, we can view an inquisitive bisimulation as a relation
which is defined exclusively on the worlds of the two models. We will
call such a relation a world-bisimulation. In order to define it, let
us first fix a way to lift a binary relation between two sets to a
relation between the corresponding powersets.

\bD
\label{bisimliftdef}
The \emph{lifting} of a relation $Y\subseteq W\times W'$ 
to information states is the relation $\overline Y \subset \wp (W)\times\wp(W')$ linking 
information states $s$ and $s'$ iff
\bre
\item[--]
 for all $w\in s$ there is a $w'\in s'\text{ s.t.\ }w Yw'$
\item[--] 
for all $w'\in s'$ there is a $w\in s\text{ s.t.\ }w Yw'$
\ere
\eD

\bD[world-bisimulation]~\\
Let $M=(W,\Sigma,V)$ and $M'=(W',\Sigma',V')$ be two inquisitive modal models. A non-empty relation $Y\subseteq W\times W'$ is called a \emph{world-bisimulation} in case  the following constraints are satisfied whenever $wYw'$:
\begin{itemize}
\item atom equivalence: 
\begin{itemize}
\item
  $\forall p\in\P$: $w\in V(p)\iff w'\in V'(p)$
\end{itemize}
\item back\&forth: 
\begin{itemize}
\item for all $s\in \Sigma(w)$ there is  $s'\in\Sigma'(w')$ s.t.\ $s\overline Ys'$
\item for all $s'\in \Sigma'(w')$ there is  $s\in\Sigma(w)$ s.t.\ $s\overline Ys'$
\end{itemize}
\end{itemize}
\eD

Bisimulations and world-bisimulations are tightly connected, 
as the following proposition brings out. The straightforward proof is omitted.

\bP 
If $Z$ is a bisimulation between two models $\M$ and $\M'$, then its
restriction to worlds, $Z^{\w} :=Z\cap(W\times W')$, 
is a world-bisimulation. 
Conversely, if $Y$ is a world-bisimulation, then $Y\cup\overline Y$ is a bisimulation. 
\eP

If $Z$ is a bisimulation, then $Z$ is included in $Z^{\w}\cup\overline{Z^{\w}}$, 
but not necessarily identical to it. 
Thus, a bisimulation is not uniquely determined by its restriction to worlds. 
Rather, given a world-bisimulation $Y$, the bisimulation
$Y\cup\overline Y$ is the largest among the bisimulations $Z$ with 
$Z^{\w}=Y$.

\bC\label{cor:bisimilarity and world-bisimulation}
Two world-pointed models $\M,w$ and $\M',w'$ are 
bisimilar iff there is a world-bisimulation $Y$ such that $wYw'$. 
Two state-pointed models $\M,s$ and $\M',s'$ are bisimilar iff there
is a world-bisimulation $Y$ such that $s\overline Ys'$.
\eC

We now turn to the finite levels of bisimilarity.

\bD 
\label{ellbisimdef}
Let $\M$ and $\M'$ be two inquisitive modal models. A back-and-forth
system of height $n$ is a family $(Z_i)_{i\le n}$ of non-empty relations 
$Z_i\subseteq W\times W'\cup\wp(W)\times\wp(W')$ satisfying the 
following constraints for each $i\le n$:
\begin{itemize}
\item atom equivalence: if $wZ_iw'$ then for all $p\in\P$, $w\in V(p)\iff w'\in V'(p)$
\item state-to-world back\&forth: if $sZ_is'$ then 
\begin{itemize}
\item for all $w\in s$ there is some $w\in s'$ s.t.\ $wZ_iw'$
\item for all $w'\in s'$ there is some $w\in s$ s.t.\ $wZ_iw'$
\end{itemize}
\item world-to-state back\&forth: if $i>0$ and $wZ_i{w'}$ then 
\begin{itemize}
\item for all $s\in \Sigma(w)$ there is some $s'\in\Sigma'(w')$ s.t.\ $sZ_{i-1}s'$
\item for all $s'\in \Sigma'(w')$ there is some $s\in\Sigma(w)$ s.t.\ $sZ_{i-1}s'$
\end{itemize}
\end{itemize}
\eD

It is straightforward to check that $n$-bisimilarity can be
characterised in terms of back\&forth systems as follows.

\bP 
Let $\M,x$ and $\M',x'$ be two world- or state-pointed models. 
$\M,x\,\simn\,\M',x'$ iff there exists a back\&forth system 
$(Z_i)_{i\leq n }$ such that $xZ_n x'$. 
\eP

Analogously to what we did for full bisimilarity, it is also possible to give a purely world-based notion of back\&forth-system of height $n$ as a family of relations $(Y_i)_{i\le n}\subseteq W\times W'$. As expected, $n$-bisimilarity can then be characterised in terms of the existence of such a system, in a way analogous to the one given by Corollary~\ref{cor:bisimilarity and world-bisimulation}. We leave the details to the reader.

\section{An Ehrenfeucht--Fra\"\i ss\'e theorem} 
\label{sec:ef}

The crucial r\^ole of these notions of equivalence for the model
theory of inquisitive modal logic is brought out 
in a corresponding Ehrenfeucht--Fra\"\i ss\'e theorem.

Using the standard notion of the modal depth of a formula, 
$\inqbm_n$ denotes the class of $\inqbm$-formulae of depth up to $n$. 
It is easy to see that the
semantics of any formula in $\inqbm_n$
is preserved under $n$-bisimilarity; as a consequence, all of inquisitive modal logic is
preserved under full bisimilarity. The following analogue
of the classical Ehrenfeucht--Fra\"\i ss\'e theorem shows that, for
finite sets \P\ of atomic propositions, 
{$n$-bisimilarity} coincides with logical 
indistinguishability in $\inqbm_n$, which we denote as 
$\equiv^n_{\inqbm}$: 
\[
\M,s \equiv^n_{\inqbm} \M',s' 
\;:\Longleftrightarrow\; \left\{\barr{@{\;}r@{}}
\M,s \models \phi\;\;\Leftrightarrow\;\;
\M',s' \models \phi
\\
\mbox{for all } \phi \in \inqbm_n.\!\!\nt \earr\right.
\]

\bT[Ehrenfeucht--Fra\'\i ss\'e theorem for \inqbm]\label{EFthm}~\\
\emph{Over any finite set of atomic propositions \P, 
for any $n \in \N$ and inquisitive state-pointed modal models $\M,s$ and $\M',s'$:}
 \bre
 \item $\M,s\;\simn\; \M',s' \iff
\M,s \equiv^n_{\inqbm}\M',s'$
\item $\M,w\;\simn\; \M',w' \iff
\M,w \equiv^n_{\inqbm}\M',w'$
\ere
\eT

Notice that item (ii) of the theorem follows from item (i) by taking
$s$ and $s'$ to be singleton states.
As usual, the crucial
implication of the theorem, from right to left, follows from the existence of
\emph{characteristic formulae} for $\sim^n$-classes of
pointed models---and it is here that the
finiteness of \P\ is crucial. 
Notice, however, that while we can expect a formula to uniquely
characterise the $\sim^n$ class of a world, we cannot expect a formula
to uniquely define the $\sim^n$-class of an information state, for
this would conflict with the persistency property of the logic 
(Proposition~\ref{prop:persistency}): if a formula is supported at $\M,s$, 
it must also be supported at $\M,s'$ for all $s'\subseteq s$ even when 
$\M,s'\not\sim^n\M,s$. However, the next proposition shows that 
$\inqbm_n$-formulae characterise the $\sim_n$-class of an information 
state up to persistency.

\bP[characteristic formulae for $\simn$-classes]
\label{prop:characteristic formulae}~\\  
Let $\M,w$ be a world-pointed model and $\M,s$ a state-pointed model 
over a finite set of atomic propositions $\P$. There are
$\inqml$-formulae $\chi^n_{\M,w}$ and $\chi_{\M,s}^n$ of modal depth $n$ such that: 
\bre
\item $\M',w'\models\chi^n_{\M,w}\iff \M',w'\,\simn\, \M,w$
\item $\M',s'\models\chi^n_{\M,s} \iff \M',s'\,\simn\, \M,t$ for some  $t\subseteq s$ 
\ere
\eP

These results can be extended straightforwardly to a multi-modal 
inquisitive setting with a finite set $\mathcal A$ of agents.

\prf
By simultaneous induction on $n$, we define formulae 
$\chi^n_{\M,w}$ and $\chi^n_{\M,s}$ together with auxiliary
formulae $\chi^n_{\M,\Pi}$ for all worlds $w$, information states $s$ and inquisitive
states $\Pi$ over $\M$. Given two inquisitive states $\Pi$ and $\Pi'$
in models $\M$ 
and $\M'$, 
we write $\M,\Pi \sim^n \M',\Pi'$ if every  
state $s\in \Pi$ is $n$-bisimilar to some 
state $s'\in\Pi'$, and vice versa. 
Dropping reference to the fixed $\M$, 
we let:
\begin{eqnarray*}\label{def:characteristic formulae}
\chi^0_w&=&\bigwedge
\{p\colon w\in V(p) \}
\land\bigwedge
\{\neg p\colon w\not\in V(p) \}
\\
\chi^n_s&=&\bigvee\{\chi^n_w\colon w\in s\}
\\[.1cm]
\chi^n_\Pi&=&\Lori\{\chi^n_s \colon s\in\Pi\}
\\
\hnt\chi^{n+1}_w&=&\chi^n_w\land\ibox\chi^n_{\Sigma(w)}\!\land\bigwedge\{\neg\!\ibox\!\chi^n_\Pi\colon
                            \Pi\subseteq\Sigma(w),\, 
\Pi\not\!\!{\simn}\Sigma(w)\}
\end{eqnarray*}

As a special case, we have $\chi^n_{\emptyset}=\bot$ 
(as $\bigvee \emptyset \equiv \bot$).
These formulae are of the required modal depth; the conjunctions and
 disjunctions in the definition are well defined since, for a
 given $n$, there are only finitely many distinct formulae of the
 form $\chi^n_w$, and analogously for $\chi^n_s$ or $\chi^n_\Pi$
 (indeed, it is easy to check that, for finite $\P$,
$\inqml_n$ is finite up to logical equivalence).
Note that, by Proposition~\ref{prop:truth-conditionality}, 
the formulae $\chi_w^n$ and $\chi_s^n$ are truth-conditional for all $n$.

We show the following: 
\begin{enumerate}
\item $\M',w'\models\chi^n_{\M,w}\iff \M',w'\,\simn\, \M,w$
\item $\M',s'\models\chi^n_{\M,s} \iff \M',s'\,\simn\, \M,t$ for some  $t\subseteq s$ 
\item $\M',s'\models\chi^n_{\M,\Pi}\iff \M',s'\,\simn\, \M,s$ for some $s\in\Pi$
\end{enumerate}

We first show that, for each individual $n$, 
(1)~$\Rightarrow$~(2)~$\Rightarrow$~(3). 
The three claims are then established together by induction on~$n$.  

For (1)~$\Rightarrow$~(2), suppose $\M',s'\models\chi_{\M,s}^n$. 
By persistency (Proposition~\ref{prop:persistency}), $\chi_{\M,s}^n$ 
is true at each $w'\in s'$; that is, for all $w'\in s'$ we have 
$\M',w'\models\bigvee\{\chi_{\M,w}^n \colon w\in s\}$. 
Since connectives have the standard behaviour in terms of
truth-conditions, this means that for any $w'\in s'$ we have
$\M',w'\models\chi_{\M,w}^n$ 
for some $w\in s$. By (1), this means that any world in $s'$ is
$n$-bisimilar to some world in~$s$. Letting $t$ be the set of worlds
in $s$ that are $n$-bisimilar to some world in $s'$, we have $t\subseteq s$ 
and  $\M',s'\sim_n \M,t$. 
Conversely, suppose $\M',s'\sim_n \M,t$ for some $t\subseteq s$. 
Then every $w'\in s'$ is $n$-bisimilar to some $w\in s$. 
By (1), this means that $\M',w'\models\chi_{\M,w}^n$, which implies 
$\M',w'\models\chi_{\M,s}^n$. Since this holds for any $w'\in s'$, 
and since $\chi_{\M,s}^n$ is a truth-conditional formula 
(by Proposition~\ref{prop:truth-conditionality}), 
it follows that $\M',s'\models\chi_{\M,s}^n$. 

For (2)~$\Rightarrow$~(3), suppose $\M',s'\models\chi_{\M,\Pi}^n$. This implies 
$\M',s'\models\chi_{\M,s}^n$ for some $s\in\Pi$. By claim (2) we have 
$\M',s'\sim_n \M,t$ for some $t\subseteq s$. Since $\Pi$ is downward
closed, $t\in\Pi$. Conversely, suppose $\M',s'\sim_n \M,t$ for some $t\in\Pi$. 
By (2), $\M',s'\models\chi_{\M,t}^n$, and since $t\in\Pi$, also 
$\M',s'\models\chi_{\M,\Pi}^n$.

We can now show (1) (and thus~(2) and~(3)) for all
$n\in\mathbb{N}$ by induction. 
The claim $\M',w'\models\chi_{\M,w}^0\Leftrightarrow\M',w'\!\sim_0\! \M,w $ 
follows immediately from the definition of $\chi_{\M,w}^0$. Now assume 
that claim (1), and thus also claims (2) and (3), hold for~$n$, and
let us consider the claim for $n+1$.

For the right-to-left direction,  suppose $\M',w'\sim_{n+1} \M,w$. We
want to show that $\M',w'\models\chi_{\M,w}^{n+1}$. This amounts to
showing that: 
\bre
\item
$\M',w'\models\chi_{\M,w}^{n}$;
\item
$\M',w'\models\ibox\chi_{\M,\Sigma(w)}^{n}$;
\item
$\M',w'\models\neg\!\ibox\chi^{n}_{\M,\Pi}$ when
$\Pi\subseteq\Sigma(w)$ and $\Pi\not\sim_n\Sigma(w)$. 
\ere

For~(i): $\M',w'\sim_{n+1} \M,w$ implies $\M',w'\sim_n \M,w$, so by the 
induction hypothesis $\M',w'\models\chi^n_{\M,w}$.

For~(ii) take $s'\in\Sigma'(w')$. Since $\M',w'\sim_{n+1} \M,w$ we must have 
$\M',s'\sim_n \M,s$ for some  $s\in\Sigma(w)$. By the induction 
hypothesis, $\M',s'\models\chi_{\M,\Sigma(w)}^n$. This holds for all 
$s'\in\Sigma'(w')$, and so $\M',w'\models\ibox\chi^n_{\M,\Sigma(w)}$.

For~(iii) suppose for a contradiction that for some
  $\Pi\subseteq\Sigma(w)$,  $\Pi\not\sim_n\Sigma(w)$  and
  $\M',w'\models\ibox\chi^n_{\M,\Pi}$. This means that every
  $s'\in\Sigma'(w')$ supports $\chi^n_{\M,\Pi}$ and thus, by our
  induction hypothesis, is $n$-bisimilar to some $s\in\Pi$. Since
  $\Pi\subseteq\Sigma(w)$ and $\Pi\not\sim_n\Sigma(w)$, there must be
  a state $t\in\Sigma(w)$ which is not $n$-bisimilar to any
  $s\in\Pi$. But since any state $s'\in\Sigma'(w')$ is $n$-bisimilar
  to some $s\in\Pi$, this means that $t$ is not $n$-bisimilar to any
  $s'\in \Sigma'(w')$. Since $t\in\Sigma(w)$, this contradicts the
  assumption that $\M',w'\sim_{n+1}\M,w$.

This establishes the right-to-left direction of claim~(1). For the converse, 
suppose $\M',w'\models\chi^{n+1}_{\M,w}$. To prove $\M',w'\sim_{n+1}\M,w$, 
we must show that: 
\bre
\item
$w'$ and $w$ coincide on atomic formulae;
\item
any $s'\in\Sigma'(w')$ is $n$-bisimilar to some $s\in\Sigma(w)$;
\item
any $s\in\Sigma(w)$ is $n$-bisimilar to some $s'\in\Sigma'(w')$.
\ere

For~(i): Since $\chi_{\M,w}^n$ is a conjunct of $\chi^{n+1}_{\M,w}$, by the induction 
hypothesis we have $\M',w'\sim_n \M,w$, which implies that $w$ and $w'$ 
satisfy the same atomic formulae.

For~(ii): Since $\ibox\chi^n_{\M,\Sigma(w)}$ is a conjunct of 
$\chi^{n+1}_{\M,w}$, $\M',w'\models\ibox\chi^n_{\M,\Sigma(w)}$. 
This implies that any $s'\in\Sigma'(w')$ supports
$\chi^n_{\M,\Sigma(w)}$. 
By induction hypothesis, this means that any $s'\in\Sigma'(w')$ is 
$n$-bisimilar to some $s\in\Sigma(w)$.

In preparation for~(iii), consider the set $\Pi$ of states in
$\Sigma(w)$ that are $n$-bisimilar to some $s'\in\Sigma'(w')$. 
We have already seen that any $s'$ is $n$-bisimilar to some state $s\in\Sigma(w)$, 
which must then be in $\Pi$ by definition. By induction hypothesis, 
the fact that $s'$ is $n$-bisimilar to some state in $\Pi$ implies 
$\M',s'\models\chi^n_{\M,\Pi}$. As this holds for every 
$s'\in\Sigma'(w')$, we have $\M',w'\models\ibox\chi^n_{\M,\Pi}$. 

Now suppose towards a contradiction that, contrary to~(iii), 
some $s\in\Sigma(w)$ were not $n$-bisimilar to any state in $\Sigma'(w')$.
Then $s$ could not be $n$-bisimilar to any state in $\Pi$ either. This
implies that $\Pi\not\sim_n\Sigma(w)$ so that $\neg\ibox\chi^n_{\M,\Pi}$ 
would be a conjunct of $\chi^{n+1}_{\M,w}$. 
Then, since $\M',w'\models\chi^{n+1}_{\M,w}$, we should have
$\M',w'\models\neg\ibox\chi^n_{\M,\Pi}$, contrary to what we found above. 
This completes the proof. 
\eprf

It is now easy to prove the 
non-trivial direction of Theorem \ref{EFthm}. 

\prf[Proof of Theorem~\ref{EFthm}]
We focus on the left-to-right direction in claim~(i) of 
Theorem~\ref{EFthm}: the converse follows
from the observation that \inqml-formulae of depth up to $n$ are
invariant under $n$-bisimilarity, and claim~(ii) follows
from~(i) by specialisation to singleton states. So suppose
$\M,s\not\sim_n \M',s'$: then either of the states $s$ and $s'$ is not
$n$-bisimilar to any subset of the other. Without loss of generality,
suppose it is~$s'$. By the property of the formula $\chi^n_{\M,s}$ we
have $\M,s\models\chi^n_{M,s}$ but
$\M',s'\not\models\chi^n_{\M,s}$. Since the  modal depth of
$\chi^n_{\M,s}$ is $n$, this shows that
$\M,s\not\equiv^n_{\inqbm}\M',s'$. 
\eprf

As a corollary of Theorem~\ref{EFthm}, we have the following
characterisation of properties definable in \inqml.

\bC 
\label{EFcorrworldpointed}\label{EFcorrstatepointed} 
A property of world-pointed models (resp., an inquisitive property of state-pointed models) over a finite set $\P$ of atomic propositions is definable in \inqbm\ if and
only if it is closed under $\simn$ for some $n\in\N$.
\eC

\prf
If a property $\KK$ of pointed models is defined by a formula $\phi$
of depth $n$, then, since $\phi$ is invariant under $n$-bisimilarity,
$\KK$ is closed under $\sim^n$. 

Conversely, suppose that $\KK$ is a property of world-pointed models closed under $\sim^n$. Using Proposition~\ref{prop:characteristic formulae} 
it is easy to show that $\KK$ is defined by the formula 
$\chi_n^{\KK}:=\bigvee\{\chi^n_{\M,w} \colon {\<\M,w\>\in{\KK}}\}$. 
Notice that the disjunction is well defined, since for a given $n$
there are only finitely many distinct formulae of the form $\chi^n_{\M,w}$. 

Similarly, if $\KK$ is an inquisitive property of state-pointed models
closed 
under $\sim^n$, it follows from 
Proposition~\ref{prop:characteristic formulae}  that $\KK$ is defined 
by the inquisitive disjunction 
$\chi_n^{\KK}=\Lori\{\chi^n_{\M,s}\colon {\<\M,s\>\in{\KK}}\}$.
\eprf

\bR\label{rem:kbox} Notice that the construction of characteristic
formulae does not use the modality $\kbox$. This implies that $\kbox$ 
can be eliminated from the language of \inqml\ without loss of expressive power. 
This was proved in a more direct way in~\cite{IvanoDiss}, where it is shown 
that a formula $\kbox\phi$ can always be turned into an equivalent 
$\kbox$-free formula. However, this translation is not schematic, i.e., 
there is no $\kbox$-free formula $\psi(p)$ such that for every 
$\phi$, $\kbox\phi\equiv\psi(\phi)$.
\eR

\section{Interlude: \inqml\  and neighbourhood semantics}
\label{sec:neighbourhood}

In neighbourhood semantics for modal logic (see~\cite{Pacuit:17} for a
recent overview), modal formulae are interpreted with respect to
\emph{neighbourhood models}, which are defined as triples 
$\M=(W,\Sigma,V)$ 
where $W$ is a set of worlds, $V:\P\to\wp(W)$ is a propositional
valuation, and 
$\Sigma:W\to\wp\wp(W)$, 
called a \emph{neighbourhood map}, is a function which assigns to each
world a set of information states. The standard language of modal
logic is interpreted on such models by means of the standard
truth-conditional clauses for connectives, and the following clause
for modalities:
$$\M,w\models_{\ssc \mathrm{nhd}} \kbox\phi\iff |\phi|_M\in\Sigma(w)$$
where $|\phi|_M$ is the set of worlds in $\M$ where $\phi$ is true. A class of neighbourhood models which is particularly well-studied is that of \emph{monotonic} neighbourhood models \cite{Hansen:03}, which are characterised by the fact that, for all worlds $w$, the set $\Sigma(w)$ is upward-closed, i.e., closed under supersets.

Clearly, an inquisitive modal model is a special case of neighbourhood model: it is a neighbourhood model such that $\Sigma(w)$ is non-empty and downward closed, i.e., closed under subsets. That is, inquisitive modal models are neighbourhood models which have exactly the opposite monotonicity property than monotonic neighbourhood models have. 

In spite of this similarity in models, however, there are big
differences between \inqml\ and neighbourhood semantics, in terms of 
the logics that arise from these approaches, their expressive power 
and the induced notions of equivalence.

These differences arise from the way in which the neighbourhood
function is used to interpret modal formulae. In neighbourhood semantics, to 
interpret $\kbox\phi$ we check whether the interpretation of $\phi$ 
is a neighbourhood. The clause for the main modality of \inqml, $\ibox$, 
is very different: just as in Kripke semantics, we have to check whether $\phi$ 
holds in all successors of the given world---only, these successors are now 
information states rather than worlds.  
As a consequence of this, whereas neighbourhood semantics gives rise 
to non-normal modal logics, the logic of the $\ibox$ modality in \inqml\ is normal: 
it validates the K axiom, as well as distribution over conjunction and
the necessitation rule.%
\footnote{In this discussion, we have set aside the modality $\kbox$ 
of \inqml\ for simplicity since, as remarked above, this modality is 
definable from $\ibox$ and the connectives. However, as shown 
in~\cite{IvanoDiss}, $\kbox$ is also a normal modality in \inqml.}

Besides giving rise to very different modal logics, \inqml\ and neighbourhood 
semantics are also different, and in fact incomparable, in terms of their 
expressive power. To see that neighbourhood semantics can draw 
distinctions that \inqml\ cannot draw, consider the formula $\kbox\top$. 
In neighbourhood semantics, this expresses the property of having the 
whole universe as a neighbourhood: 
\[
\M,w \models_{\ssc\mathrm{nhd}} \kbox\top\iff W\in \Sigma(w)
\]
This property is clearly not invariant under inquisitive bisimulations
(indeed, it 
is not preserved under disjoint unions!). Thus, it is not expressible in \inqml.

To see that \inqml\ can also draw distinctions that neighbourhood semantics 
cannot draw, consider the formula $\ibox{?p}$. This formula expresses the 
fact that at every neighbourhood of the evaluation world, the truth-value 
of $p$ is constant.
$$\M,w\models\ibox{?p}\iff\forall s\in \Sigma(w): 
s\subseteq|p|_\M\text{ or }s\subseteq|\neg p|_M$$

We claim that this property is not expressible in neighbourhood semantics. 
To see this, consider two models $\M_1$ and $\M_2$ with the same 
universe $W=\{v,u,u'\}$ and the same valuation $V(p)=\{v\}$. 
The two models differ in their neighbourhood map, which are both
constant, with values
\[
\{\{v\},\{u\}\}^\downarrow \; \mbox{ for $\Sigma_1$ in $\M_1$ \quad versus } \quad
\{\{v,u\}\}^\downarrow \; \mbox{ for $\Sigma_2$ in  $\M_2$.}
\]

Given any $w\in W$, we have  $\M_1,w\models\ibox{?p}$ but 
$\M_2,w\not\models\ibox{?p}$. However, the set $\{v,u\}$ is not the truth-set of any formula
in neighbourhood semantics: the reason is that $u$ and $u'$ are indistinguishable, and so a truth-set always contains either both of them, or neither. 
Using this fact, we can show by induction that 
$\M_1,w\models_{\ssc\mathrm{nhd}} \phi\iff\M_2,w\models_{\ssc\mathrm{nhd}} \phi$ for all formulae $\phi$. 
Hence, the property expressed by $\ibox{?p}$ is not 
expressed by any formula in neighbourhood semantics.

Clearly, since \inqml\ and neighbourhood semantics are sensitive 
to different features of a model, the appropriate notion of bisimilarity is also 
different in these two contexts. For instance, consider again the above models 
$\M_1$ and $\M_2$: according to the notion of bisimilarity $\sim^{\text{N}}$ 
appropriate for neighbourhood semantics \cite{Hansen:09}, 
the relation $R=\{\<v,v\>,\<u,u\>,\<u,u'\>,\<u',u\>,\<u',u'\>\}$ is a bisimulation. 
This implies that $\M_1,v\sim^{\text{N}}\M_2,v$. By contrast, a single round
of the inquisitive bisimulation game suffices to show that 
$\M_1,v\not\sim\M_2,v$ in our setting. 

Conversely, under our notion of bisimulation, a point $w$ in a model $\M$ 
is always fully bisimilar to its copy in the disjoint union $\M\uplus\M'$. 
Clearly, the same cannot be true in neighbourhood semantics, 
given that in this semantics modal formulae 
are not in general preserved under disjoint unions.

A notion of bisimulation which is much closer to the one we study here is found in
the literature on monotonic neighbourhood models \cite{Hansen:03}. 
In terms of the bisimulation game, the difference between the two notions 
can be described as follows.
Starting from a world-position $\<w,w'\>$, Player \PlayerI\ picks a state $s$ in either 
$\Sigma(w)$ or $\Sigma'(w')$; Player \PlayerII\ responds with a state $s'$ on 
the opposite side. At this point, the two games come apart: in our
version of the game, \PlayerI\ can choose a world from either $s$ or $s'$, while in
the version given in~\cite{Hansen:03}, \PlayerI\ is required to pick a world from $s'$. 
Imposing such a restriction in our setting would trivialise the game, 
providing \PlayerII\ with a universal winning strategy: always pick
$s'=\emptyset$.

Interestingly, however, one can show that due to the downward-closure of $\Sigma(w)$, in our setting it would not make a difference (in terms of the resulting notion of bisimilarity) if Player \PlayerI\ were required to pick a world from the state $s$ that he himself selected in the world-to-state phase. Thus, we could equivalently have presented the game in a form which is the mirror image of the game used in monotonic neighbourhood frames. Clearly, this symmetry reflects the opposite monotonicity constraints that these two logics place on the neighbourhood map. 

\section{Relational inquisitive models}
\label{relinqmodsec}

In the remainder of this paper
we compare the expressive power of
inquisitive modal logic with that of first-order logic. This 
is not quite as straightforward as for ordinary modal logic. 
A standard Kripke model can be identified naturally with a relational structure
with a binary accessibility relation $R$ 
and a unary predicate for the interpretation of each atomic 
proposition $p \in \P$. 
By contrast, an inquisitive modal model also needs to encode the
inquisitive state map $\Sigma:W\to\wp\wp(W)$. This map can be identified with a binary relation 
$E \subseteq W\times \wp(W)$. In order to view this as part of a
relational structure, however, we need to adopt a two-sorted perspective, 
and view $W$ and $\wp(W)$ as domains of two distinct sorts. We thus turn to two-sorted structures. In order to capture the fact that the second sort contains sets of elements of the first sort, our relational structures include a relation $\epsilon$ between these sorts, which simulates set-theoretic membership. 
This leads  to the following notion.

\subsection{Relational inquisitive modal models}

\bD[relational models]
\label{def:relational models}~\\
A \emph{relational inquisitive modal model} over atomic 
propositions $\P= \{ p_i \colon i\in I \}$ is a relational structure
\[
\str{M}=(W,S,E,\ee,(P_i)_{i\in I})
\]
where $W,S$ are non-empty sets related by $E,\ee\subseteq W\times S$,
and $P_i\subseteq W$ for $i\in I$. 

With $s\in S$ we associate the set 
$\underline{s} :=\{w\in W \colon (w,s) \in \ee \} \subset W$ 
and require the following conditions, which enforce resemblance 
with inquisitive modal models:
\begin{itemize}
\item extensionality: 
for $s,s' \in S$, $\underline{s}=\underline{s}'$ implies $s=s'$.
\item local powerset: if $s\in S$ and $t\subseteq\underline{s}$, there is an $s'\in S$ such that 
$\underline{s}'=t$.
\item non-emptiness: $E[w]\neq\emptyset$ for all $w \in W$.
\item downward closure:
for $s,s' \in S$ with 
$\underline{s'}\subseteq\underline{s}$, 
$s\in E[w]$ implies $s'\in E[w]$.
\end{itemize}
Multi-modal variants are analogously defined, with a relation 
$E_a \subset W \times S$ to encode the inquisitive assignments
$\Sigma_a$ for agent $a \in \A$.

By a \emph{world (resp.\ state)-pointed} relational model we mean a pair $\<\str{M},x\>$ where $x$ is an element in the first (resp.\ second) sort of $\str{M}$.
\eD

By extensionality, the second sort $S$ of such a relational model can 
always be identified with a domain of sets over the first sort,
namely, 
$\{\underline{s} \colon s\in S\}
\subset \wp(W)$. 
In the following, we will assume this identification and view a 
relational model as a structure $\str{M} =(W,S,E,\in,(P_i))$ where 
$S\subseteq\wp(W)$ and $\in$ is the actual membership relation.
We shall therefore also specify relational
models by just $\str{M}=(W, S, E,(P_i))$ when the fact that $S \subset
\wp(W)$ and the natural interpretation of $\ee$ are understood.
Notice that, given this identification, the downward closure condition
can be stated more simply as: if $s\in E[w]$ and $t\subseteq s$, then $t\in E[w]$.

Notice that a relational model $\str{M}$ induces a corresponding Kripke structure $\str{K}(\str{M})=(W, R,(P_i)_{i\in I})$,
where $R\subseteq W\times W$ is the relation defined as follows:
\[
wRw'\iff \text{for some }s\in S:wEs\text{ and }w'\in s,
\] 
so that 
$R[w]:=\{w' \colon wRw'\}=\bigcup E[w]$ as the natural 
relational encoding of the map $\sigma \colon w \mapsto \bigcup\Sigma(w)$.

In addition to the above conditions, we might
impose other constraints on a relational model $\str{M}$: in particular, we may require
$S$ to be the full powerset of $W$, or to resemble the powerset 
from the local perspective of each world~$w\in W$.

\bD
\label{relclassesdef}
A relational model $\str{M}=(W,S,E,(P_i))$ is 
\bre
\item[--] \emph{full} if $S=\wp(W)$;
\item[--] \emph{locally full} if $\wp(R[w]) \subseteq S$ for all $w\in W$.
\ere
\eD

Note that, by local powerset, the condition of 
local fullness is equivalent to the condition that 
the information states $R[w] = \bigcup E[w]$ 
are represented in $S$ for all $w \in W$.

\subsection{Relational encoding of inquisitive modal models}
\label{relencsec}
 
The connection between inquisitive modal models and their relational
counterparts is not one-to-one. In one direction, a relational model
$\str{M}=(W,S,E,(P_i))$ uniquely determines an inquisitive modal
model 
\[
\str{M}^\ast=(W,\Sigma,V) \;\;\mbox{ where } \Sigma(w)=E[w], V(p_i)=P_i.
\] 

Notice that the non-emptiness and downward closure
conditions on $E$ guarantee that $\str{M}^\ast$ is indeed an inquisitive
modal model. 
Since the passage from $\str{M}$ to $\str{M}^\ast$ obliterates information about
the second sort $S$, there are in general many different relational
models $\str M$ that determine the same inquisitive modal model
$\M$. That is, a given inquisitive modal model may have different
relational counterparts. Let us call such counterparts the
\emph{relational encodings} of $\M$.

\bD 
\label{relencodedef}
A \emph{relational encoding} of an inquisitive modal model $\M$ is
a relational model $\str M$ with $\str M^\ast=\M$.  \eD

Clearly, two relational counterparts of $\M$ must coincide in terms of
$W$, $E$ and the $P_i$. 
But this leaves quite some choice with respect to the richness 
of the second sort. The following isolates some immediate choices.

\bD[relational encodings]
\label{relencdef}~\\
Given an inquisitive modal model $\M=(W,\Sigma,V)$, we define three 
relational encodings $\str{M}^{\ssc[\cdots ]}(\M)$ 
of $\M$, each based on~$W$, and with 
$wEs\Leftrightarrow s\in\Sigma(w)$, $w\,\ee\, s\Leftrightarrow w\in s$ 
and $P_i  = V(p_i)$. 
The encodings differ in the second sort domain~$S$:
\begin{itemize}
\item for $\MM^\rel(\M)$, the minimal encoding of $\M$: 
\\
$S:=\mathrm{image}(\Sigma)$;
\item for $\MM^\lf(\M)$, the minimal locally full  encoding of $\M$: 
\\
$S:=\{s\subseteq\sigma(w)\colon
  w\in W\}$;
\item for $\MM^\full(\M)$, the unique full encoding of $\M$: \\ $S:=\wp(W)$.
\end{itemize}
To encode state-pointed models $\M,s$ we augment the corresponding 
$S$ by $\wp(s)$. 
These definitions generalise in a natural way to the multi-modal case.
\eD

\subsection{Relational models and \inqml}
\label{standardtransec}

\noindent The notions of ($n$-)bisimilarity defined 
in Section~\ref{sec:bisimulations} naturally lifts to 
relational models as follows.

\bD[(n-)bisimilarity for relational models]~\\
Two state- or world-pointed relational models 
$\str{M}_1,x_1$ and $\str{M}_2,x_2$ are 
bisimilar, $\str{M}_1,x_1\sim\str{M}_2,x_2$, if they encode
bisimilar pointed models, i.e.\ if 
$\str{M}_1^\ast,x_1\sim\str{M}_2^\ast,x_2$.
Similarly for $n$-bisimilarity: 
$\str{M}_1,x_1\sim_n\str{M}_2,x_2$ if 
$\str{M}_1^\ast,x_1\sim_n\str{M}_2^\ast,x_2$.
\eD

In particular, any two encodings of the 
same pointed inquisitive model are bisimilar.
We similarly lift the interpretation of \inqml\ to relational models.

\bD Let $\str{M},x$ be a (world- or state-)pointed relational model
and $\phi$ a formula of \inqml. We define $\str{M},x\models\phi$
to mean $\str{M}^\ast,x\models\phi$.
\eD

By a \emph{property} of world-pointed relational models we mean a class of world-pointed relational models. Similarly, by an \emph{inquisitive property} of state-pointed relational models we mean a class of state-pointed relational modes which satisfies the analogues of the conditions in Definition~\ref{inquisitive state-properties}. By a world- or state-property over a class of relational models $\CC$ we mean a class of pointed models $(\str{M},x)$ with $\str{M}\in\CC$.
We can then translate Corollary~\ref{EFcorrworldpointed} into a characterisation of the properties of 
pointed relational models that are definable in \inqml.

\bC 
\label{EFcorrworldpointedrel}\label{EFcorrstatepointedrel} 
Let $\KK$ be a property of world-pointed relational models, or an inquisitive property of state-pointed relational models, over a finite set $\P$ of atomic propositions. Then $\KK$ is definable in \inqbm\ if and
only if it is closed under $\simn$ for some $n\in\N$.
More generally, if $\KK$ is a world property or an inquisitive state property over a class $\CC$ of relational models, $\KK$ is definable in \inqbm\ over $\CC$ 
iff it is closed under $\simn$ within $\CC$ for some $n\in\N$.
\eC

\prf It follows from Theorem~\ref{EFthm} that if a class $\KK$ is defined by an \inqml\ formula $\phi$, then $\KK$ is closed under $\sim_n$ where $n$ is the modal depth of $\phi$. 

Conversely, suppose $\KK$ is a property 
of world-pointed relational models closed under $\simn$
for some $n\in\N$. Then 
$\KK^\ast:=\{(\str{M}^\ast,w)\colon (\str{M},w)\in\KK\}$ 
is also closed under $\simn$, so by Corollary~\ref{EFcorrworldpointed} 
it is definable by a formula $\phi$ of \inqml. 
It is easy to check that $\phi$ defines $\KK$
(note that $\sim^n$-closure of $\KK$ implies 
that $(\str{M},w)\in\KK\iff (\str{M}^*,w)\in\KK^*$).
If $\KK \subseteq \CC$ is closed under $\simn$ in restriction to
$\CC$, we may similarly work with 
$\KK^\ast:=\{(\str{M}^\ast,w)\colon 
\str{M},w \,\simn\, \str{M}',w' \mbox{ for some }
(\str{M}',w') \in\KK\}$.  
The reasoning for inquisitive state-pointed properties 
is exactly analogous, 
using the second part of Corollary~\ref{EFcorrworldpointed}. 
\eprf

\subsection{Relational models and first-order logic}
\label{standardtransec}

A relational inquisitive model supports a two-sorted first-order
language having two relation symbols $\E$ and $\eee$ corresponding
to the relations $E$ and $\ee$, respectively, and predicate symbols $\PP_i$ for $i\in I$. 
We use $\w,\v,\texttt{u}$ as variables for the first sort, and $\s,\t$ as
variables for the second sort.%
\footnote{We use different
  fonts to distinguish object language symbols ($\E,\w,\s,\dots$), in
  typewriter font, from the corresponding notation for semantic objects
  ($E,w,s,\dots$), in regular~font.}
Moreover, we make use of two defined binary predicates. 
The first is simply inclusion, defined in the natural way in terms of $\eee$:
\[
\s\subseteq\t := \forall\w(\eee(\w,\s)\to\eee(\w,\t)).
\] 

The second defined predicate, 
$\e(\w,\t)$, corresponds to the relation $R[w]=t$
(i.e., the relational encoding of the graph of the map $\sigma$): 
\[
\e(\w,\t):=\forall\v(\eee(\v,\t)\leftrightarrow
\exists\s(\E(\w,\s)\land\eee(\v,\s))).
\]

In terms of this language we can define a pair of standard
translations $\stw(\phi)$ of $\sts(\phi)$ of a formula, which capture
its truth conditions in a world and its
support conditions in an information state, respectively. 
Correspondingly, $\stw(\phi)$ has a single free variable $\w$ 
of the first sort while $\sts(\phi)$ has the free variable $\s$ of the
second sort. 
Of $\sts(\phi)$ we also use a substitution variant $\stt(\phi)$  
which is just like $\sts(\phi)$ except that the roles of variables
$\s$ and $\t$ are exchanged. 
The following define these standard translations by simultaneous induction:

\begin{itemize}
\item $\stw(p_i)=\texttt P_i(\w)$\\
$\sts(p_i)=\forall \w(\eee(\w,\s)\to\stw(p_i))$
\item 
$\stw(\bot)=\bot$\\
$\sts(\bot)=\forall \w(\eee(\w,\s)\to\stw(\bot)) \quad 
\bigl( \equiv \neg \exists \w\, \eee(\w,\s) \bigr)$
\item $\stw(\phi\land\psi)=\stw(\phi)\land\stw(\psi)$\\
$\sts(\phi\land\psi)=\sts(\phi)\land\sts(\psi)$
\item $\stw(\phi\lori\psi)=\stw(\phi)\lor\stw(\psi)$\\
$\sts(\phi\lori\psi)=\sts(\phi)\lor\sts(\psi)$
\item $\stw(\phi\to\psi)=\stw(\phi)\to\stw(\psi)$\\
$\sts(\phi\to\psi)=\forall\t(\t\subseteq\s\to(\stt(\phi)\to\stt(\psi)))$
\item $\stw(\ibox\phi)=\forall\s(\E(\w,\s)\to\sts(\phi))$\\
$\sts(\ibox\phi)=\forall\w(\eee(\w,\s)\to\stw(\ibox\phi))$
\item 
$\stw(\kbox\phi)=\forall\s(\e(\w,\s)\to\sts(\phi))$\\
$\sts(\kbox\phi)=\forall\w(\eee(\w,\s)\to\stw(\kbox\phi))$\\
\end{itemize}

It is straightforward to verify that the truth-conditions and support-conditions 
of $\phi$ in a model $\M$ correspond, respectively, to the satisfaction 
conditions for $\stw(\phi)$ and $\sts(\phi)$ in any locally full 
relational encoding of $\M$.

\bP Let 
$\str{M}$ be a locally full relational inquisitive model, $\phi\in\inqml$. 
For all worlds $w\in W$ and all states $s\in S$:
\bre
\item $\MM,w\models\phi\iff\MM,w\models\text{\emph{ST}}_{\emph{\texttt{w}}}(\phi)$
\item $\MM,s\models\phi \iff \MM,s\models\text{\emph{ST}}_{\emph{\texttt{s}}}(\phi)$
\ere
\eP

The assumption that $\MM$ be locally full is crucial for
this result. This is because, if a model is not locally full, then for
some $w\in W$ it could be that the state $\sigma(w)=\bigcup\Sigma(w)$
which is involved in determining the truth condition of $\kbox\phi$ is
not represented in $\MM$. If so, there will be no state $s\in S$
satisfying $\eee(w,s)$, which means that $\stw(\kbox\phi)$ will come
out as vacuously true at $w$, regardless of whether or not
$\MM^\ast,w\models\kbox\phi$.
However, even when $\MM$ is not locally full, 
preservation still holds for all $\kbox$-free formulae, as one can
easily verify.

\bP
\label{prop:correspondence kbox-free} 
Let $\str{M}$ be a relational inquisitive 
model, $\phi\in\inqml$ a $\kbox$-free formula. 
Then for all worlds $w\in W$ and all states $s\in S$:
\bre
\item $\MM,w\models\phi\iff\MM,w\models\text{\emph{ST}}_{\emph{\texttt{w}}}(\phi)$
\item $\MM,s\models\phi \iff \MM,s\models \text{\emph{ST}}_{\emph{\texttt{s}}}(\phi)$
\ere
\eP

Recall that, by Remark~\ref{rem:kbox}, any formula $\phi$ of \inqml\ is 
equivalent to some 
$\kbox$-free formula $\phi^*$. Combining this with the previous proposition, 
we have the following corollary.

\bC\label{cor:correspondence full}
For any $\phi\in\inqml$ there exist first-order
formulae $\phi_{\emph{\texttt{w}}}^\star:=\text{\emph{ST}}_{\emph{\texttt{w}}}(\phi^*)$ 
and $\phi_{\emph{\texttt{s}}}^\star:=\text{\emph{ST}}_{\emph{\texttt{s}}}(\phi^*)$
such that for any relational 
inquisitive model $\MM$, world $w \in W$ and $s \in S$:
\bre
\item $\MM,w\models\phi\iff\MM,w \models\phi_{\emph{\texttt{w}}}^\star$
\item $\MM,s\models\phi \iff \MM,s\models\phi_{\emph{\texttt{s}}}^\star$
\ere
\eC

The corollary allows us to view $\inqbm$ as a
syntactic fragment of first-order logic, $\inqbm \subset \FO$, over the class 
of all relational inquisitive models, just as standard
modal logic $\ML$ may be regarded as a fragment $\ML\subset \FO$ 
over Kripke models. 
Importantly, however, the class of relational inquisitive modal models
is not first-order definable in this framework, since the local 
powerset condition involves a second-order quantification. In other
words, we are dealing with first-order logic 
over non-elementary classes 
of intended models. In fact, first-order logic is not compact over this class, as the following example shows.

\bE
\label{compfailrelex}
There is a first-order formula $\phi(\s)$ in a single free variable $\s$  of
the second sort (information state) which over any relational
inquisitive model says of an element $s$ that 
there are no infinite $R$-paths inside $s$.
Combining this, for instance, with a formula that says that 
$R$ in restriction to $s$ defines a discrete linear ordering with a
minimal element, and formuale $\psi_n(\s)$ saying that $s$ 
comprises at least $n$ distinct worlds, we get a violation of compactess.
\eE

\prf
The induced modal accessibility relation $R$ is definable according to
\[
R (\texttt{u},\v)  \leftrightarrow \;\exists s ( E(\texttt{u},s) \wedge \eee(\v,\s)).
\]

The local power set condition implies that the
entire power set $\wp(s)$ of the designated state $s$ is
represented in the second sort of the relational
model. So the following formula faithfully emulates the standard
monadic second-order formalisation of the relevant property: 
\[
\phi(\s)\,:=\,\neg \exists\t \bigl(\t\subseteq\s\land \exists\mathtt{u}\,\eee(\mathtt{u},\t)\wedge \forall \texttt{u}
\bigl(\eee(\texttt{u},\t) 
\rightarrow 
\exists\v( \eee(\v,\texttt{t}) \wedge R(\texttt{u},\v))\bigr)
\bigr).
\]
where again $\t\subseteq\s$ abbreviates 
$\forall\v (\eee(\v,\t) \to \eee(\v,\s))$.
\eprf

We remark that all the considerations of this section admit 
 straightforward variations for the multi-modal
inquisitive setting, where models are equipped with a family
$(\Sigma_a)_{a \in \A}$ of inquisitive assignments, indexed by a set $\A$ of agents.

In an extension and variation of the above, Silke Mei{\ss}ner \cite{Meissner}
has proposed an alternative standard translation, which in some way is 
more uniform as it allows for a direct treatment of $\Box$. As outlined 
in~\cite{MeissnerOtto}, it also relaxes the constraints on relational 
encodings so as to extend the scope of the standard translation 
to an elementary class of relational structures, which in turns gives 
rise to a model-theoretic compactness proof for $\inqml$. 
While our translation could in principle be replaced by the more recent one 
from~\cite{MeissnerOtto}, adherence to our narrower classes of natural 
relational encodings of the \emph{intended} inquisitive models can be seen as 
a strength of our characterisation theorems.

\section{Bisimulation invariance}
\label{bisiminvsec}

\subsection{Bisimulation invariance as a semantic constraint} 
\label{bisimconstraintsec}

As discussed above, 
\inqml\ can be thought of as a fragment of first-order logic when
interpreted over relational models. We may think of $\sim$-invariance 
as a characteristic semantic feature of this fragment. 
The question we are interested in is: 
with respect to what classes $\CC$ of relational models can \inqml\ 
be characterised as being \emph{exactly} the $\sim$-invariant fragment 
of first-order logic? Let us first make precise what this means.

\bD 
We say that $\inqml$ is the $\sim$-invariant fragment of $\FO$ for 
world-properties with respect to a class $\CC$ of relational models, in symbols
\[
\inqbm \equiv_\CC^{\w} \FO/{\sim},
\] 
if for every property $\KK$ of state-pointed models over $\CC$, 
$\KK$ is definable in \inqml\ if and only if it is 
both definable in $\FO$ and $\sim$-invariant.

Similarly, we say that $\inqml$ is the $\sim$-invariant fragment of $\FO$ for inquisitive state-properties with respect to $\CC$, in symbols
\[
\inqbm \equiv_\CC^{\s} \FO/{\sim},
\] 
if for any inquisitive property $\KK$ of state-pointed models over
$\CC$, $\KK$ is definable in \inqml\ if and only if it is definable in 
$\FO$ and $\sim$-invariant.
\eD

\bR
\label{statetoworldrem}
For any class 
$\CC$, $\inqbm \equiv_\CC^{\s} \FO/{\sim}$ 
implies $\inqbm \equiv_\CC^{\w} \FO/{\sim}$.
\eR

\prf
With the $\sim$-invariant world property defined by 
$\phi(\w) \in \FO$ associate the property defined by 
$\phi'(\s) = \forall \w (\w \in \s \rightarrow \phi(\w))$. 
This property is inquisitive and $\sim$-invariant. 
From $\inqbm \equiv_\CC^{\s} \FO/{\sim}$ we obtain a formula 
$\psi \in \inqbm$ expressing this property. 
Then the same formula $\psi$, on the level of worlds, 
defines the world-property defined by $\phi(\w)$. 
\eprf

Our question then can be formulated succinctly as follows: over which 
classes $\CC$ do we have $\inqbm \equiv_\CC^{\s} \FO/{\sim}$ 
(and thus also $\inqbm \equiv_\CC^{\w} \FO/{\sim}$)?
Equivalently, the question is: 
over which classes $\CC$ is $\inqbm$ sufficiently expressive to 
capture \emph{all} first-order definable properties of 
world- or state-pointed relational models that are 
invariant under inquisitive bisimulation?

Section~\ref {inqmlcharsec} will establish that $\inqml$ is
expressively complete for $\sim$-invariant first-order properties
in this sense over each of the following classes of relational modles:
the class of all relational models, 
the class of all finite relational models, 
the class of all locally full models, and 
the class of all finite locally full models.
Before delving into the proof, however, we
discuss some underlying model-theoretic concerns
and limitations. In particular we stress the connection with 
the all-important classical r\^ole of first-order compactness,
as well as the r\^ole of non-classical model-theoretic techniques 
in dealing with first-order logic over 
non-elementary classes of models.

\subsection{Bisimulation invariance and compactness} 

Let $\KK$ be a property of world-pointed relational models over a class  $\CC$. 
Suppose $\KK$ is \inqml-definable by a formula $\phi$: then $\KK$ 
is both $\FO$-definable (by the standard translation $\phi^\star$) and 
$\sim$-invariant (as it is $\sim_n$ invariant, with $n$ the modal depth of $\phi$). 
Thus, one direction of the equivalence $\inqbm \equiv_\CC^{\w} \FO/{\sim}$ 
holds for any class $\CC$. By Corollary~\ref{EFcorrworldpointedrel}, 
the converse direction amounts to the claim that if $\KK$ is $\FO$-definable 
and $\sim$-invariant, then it is in fact $\sim_n$-invariant for some $n$.%
\footnote{Notice that if $\KK$ is $\FO$-definable, the defining formula contains only finitely many atoms. Thus, the property $\KK$ depends only on the restriction of a model to a finite set $\P$ of atoms, and we can use Corollary~\ref{EFcorrworldpointedrel} to conclude that if $\KK$ is $\sim_n$-invariant for some $n$, it is \inqml-definable.} 
That is, it amounts to the claim that, for every $\FO$-definable
property $\KK$ over $\CC$, $\sim$-invariance implies 
$\sim_n$-invariance for some $n$.

Analogous reasoning establishes the same connection w.r.t.\
inquisitive state-properties. We summarise these in the following,
where a first-order formula $\phi(\s)$ in a free variable of the
second sort is called \emph{inquisitive} over the class $\CC$ 
if the state-property expressed by $\phi(\s)$ over $\CC$ is an 
inquisitive property, i.e.\ is downward closed and always holds of the
empty state. 

\bO
\label{compactobs}
For any class $\CC$ of 
relational inquisitive models, the following are equivalent:
\bre
\item
$\inqbm \equiv_\CC^{\w} \FO/{\sim}$,
\item for any formula 
 $\phi(\w) \in
\FO$ in a single free variable of the first sort, $\sim$-invariance over $\CC$ implies 
 $\simn$-invariance over $\CC$ for some finite $n$. 
\ere
Similarly, the following are equivalent:
\bre
\item
$\inqbm \equiv_\CC^{\s} \FO/{\sim}$,
\item for any formula  $\phi(\s) \in
\FO$ in a single free variable of the second sort that
is inquisitive over the class $\CC$, $\sim$-invariance over $\CC$ implies 
$\simn$-invariance over $\CC$ for some finite $n$. 
\ere
\eO

In both contexts, condition~(ii) may be viewed as a \emph{compactness principle} for 
$\sim$-invariance of first-order properties, which is non-trivial in
the non-elementary setting of relational inquisitive models. 
Interestingly, this compactness principle for $\sim$-invariance of
$\FO$-properties of worlds fails relative to the class of \emph{full}
relational models.

\bP 
\label{failcompex} 
There is a first-order formula 
$\phi(\w)$ in a single free variable $\w$  of
the first sort that, relative to the class of full relational models, is $\sim$-invariant but not $\simn$-invariant for any $n$.
\eP

\prf\ 
Compare Example~\ref{compfailrelex} for
the following well-foundedness property:
\[
\mathcal{P}(w) :=\text{ there is no infinite }R\text{-path from }w. 
\]

On one hand $\mathcal{P}$ clearly is $\sim$ invariant but not $\simn$-invariant
for any $n$. On the other hand $\mathcal{P}$ is first-order definable over the 
class of full relational inquisitive models since,
over these models, first-order logic affords the full expressive power
of monadic second-order quantification over the first sort $W$:
first-order quantification over the second sort $S = \wp(W)$
\emph{is} quantification over subsets of the first sort. The 
formula 
\[
\phi(\w):=\neg \exists\s \bigl(
\eee (\w,\s) \wedge \forall \texttt{u} \bigl( 
\eee (\texttt{u},\s) \rightarrow 
\exists \v  (\eee(\v,\s) \wedge R(\texttt{u},\v))\bigr)
\bigr).
\]
defines the world-property $\mathcal{P}$ over any full relational 
model.
\eprf

A similar well-foundedness property can also 
be captured in first-order logic over full relational 
inquisitive epistemic models for two agents and using one basic
proposition. It suffices to 
describe analogous path properties for paths formed by a strict 
alternation of $R_a$- and $R_b$-edges on a path that alternates
between worlds where $p$ is true and where $p$ is false, for  
some atomic proposition $p$ and distinct agents  $a,b \in \A$. 
This shows that $\inqml\not\equiv_\CC^{\w} \FO/\!\!\sim$ when 
$\CC$ is the class of full relational models or full relational 
epistemic models.
Over these classes, there are $\FO$-definable, $\sim$-invariant world properties
that are not $\inqml$-definable. 
Although this is in sharp contrast with our Theorem~\ref{main1}, 
the fact that the analogue of the theorem fails over full relational models is not too surprising: over such models, $\FO$, unlike $\inqbm$, 
has access to full-fledged monadic second-order quantification.

\subsection{A non-classical route to expressive completeness} 
\label{sec:route}

In all our characterisation theorems to be treated in the following
section, we establish semantic correspondences:
\[
\qquad\qquad\inqbm \equiv_\CC^{\w} \FO/{\sim} \qquad 
\quad  \inqbm \equiv_\CC^{\s} \FO/{\sim}  
\] 
These
are assertions about equal expressive power 
between two systems presented in very different style: 
while $\inqbm$ is based on concrete syntax with clearly defined
semantics, $\FO/{\sim}$ is defined in terms of the 
semantic constraint of $\sim$-invariance. In fact, 
$\sim$-invariance is easily seen to be undecidable 
as a property of first-order formulae, so that  $\FO/{\sim}$
itself cannot be regarded as a syntactic fragment. 
As discussed in the previous subsection,
proving one of these equivalences boils down to establishing a
compactness principle relating $\sim$-invariance to $\simn$-invariance 
for some finite level $n$, in the non-classical context 
of non-elementary classes of relational inquisitive models. 

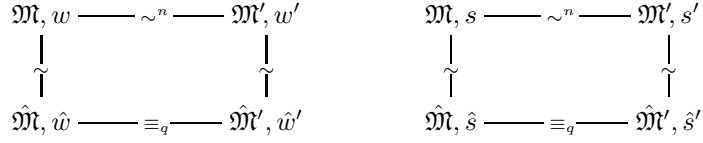
\begin{figure}
\[
\xymatrix{
\str{M},w \ar@{-}[d]|{\rule{0pt}{1ex}\sim} \ar @{-}[rr]|{\;\simn\,}
&& \str{M}'\!,w' \ar @{-}[d]|{\rule{0pt}{1ex}\sim}
\\
\hat{\str{M}},\hat{w} \ar@{-}[rr]|{\rule{0pt}{1.5ex}\;\equiv_q}
&& \hat{\str{M}}',\hat{w}'
}
\qquad\qquad
\xymatrix{
\str{M},s \ar@{-}[d]|{\rule{0pt}{1ex}\sim} \ar @{-}[rr]|{\;\simn\,}
&& \str{M}'\!,s' \ar @{-}[d]|{\rule{0pt}{1ex}\sim}
\\
\hat{\str{M}},\hat{s} \ar@{-}[rr]|{\rule{0pt}{1.5ex}\;\equiv_q}
&& \hat{\str{M}}',\hat{s}'
}
\]
\caption{Generic upgrading patterns.}
\label{genupgradefigure}
\end{figure}

For this there is a general approach that 
has been successful in a number of similar
investigations, starting from an elementary and constructive proof
in~\cite{OttoNote} of van~Benthem's characterisation theorem~\cite{Benthem83} and
its finite model theory version due to Rosen~\cite{Rosen} (for
ramifications of this method, see also~\cite{OttoAPAL04,DawarOttoAPAL09} 
and~\cite{Otto12JACM}).  
This approach involves an \emph{upgrading} of a sufficiently high 
finite level $\simn$ of bisimulation equivalence 
to a finite target level $\equiv_q$ of elementary equivalence, 
where $q$ is the quantifier rank of $\phi$. 
Concretely, and in the case of properties of worlds,
this  amounts to providing, for any world-pointed relational model
$\str{M},w$ 
a fully bisimilar pointed model $\hat{\str{M}},\hat{w}$ with the property 
that, if $\str{M},w\;\simn\;\str{M}',w'$, then  
$\hat{\str{M}},\hat{w} \equiv_q\hat{\str{M}}',\hat{w}'$.
The diagram on the left in Figure~\ref{genupgradefigure} 
shows how $\sim$-invariance of~$\phi$, together with its nature as a first-order formula 
of quantifier rank $q$, entails its $\simn$-invariance:
one chases the diagram through its lower
rung to check that, for $\phi$ that is preserved under both $\sim$ and
 $\equiv_q$, we have 
$\str{M},w \models \phi$ iff  $\str{M}',w' \models \phi$. 

The reasoning for inquisitive properties of information states is 
analogous, using a corresponding upgrading for state-pointed models
(cf.\ the right hand side in Figure~\ref{genupgradefigure}). 
At the technical level we shall mostly restrict the explicit
discussion to the more familiar world-pointed scenario, 
and only mention the necessary variations for the state-pointed case 
where relevant.

Any upgrading of the kind we just discussed 
involves an interesting tension between the very distinct 
levels of expressiveness of $\inqbm$-formulae and
$\FO$-formulae. While the latter can, for instance, distinguish 
worlds w.r.t.\ finite branching degrees of the accessibility relation $R$
or w.r.t.\ short cycles that $R$ may form in the vicinity of a world,
no $\sim$-invariant logic can. The challenge is to overcome this
discrepancy in bisimilar companion structures, using the  
malleability up to $\sim$ of relational inquisitive models 
(within the respective class $\CC$!)---and, for instance, to boost all
multiplicities and lengths of all cycles beyond what can be distinguished
in $\FO_q$ ($\FO$ up to quantifier rank $q$). 

We show how to achieve the required upgradings for various
classes $\CC$ of models in the next section, and thus
establish our characterisation theorems.

We use a variation of the upgrading technique from~\cite{OttoNote} to
instantiate the above idea. In effect we shall deviate slightly from
the generic picture in Figure~\ref{genupgradefigure} 
by interleaving $\sim$-preserving pre-processing steps and 
$\equiv_q$-preserving steps as shown in Figure~\ref{iqbmfigure}.
The upgrading itself is based on an inquisitive analogue of partial
tree unfoldings, combined with locality arguments for 
first-order Ehrenfeucht--Fra\"\i ss\'e games. We start 
with two technical remarks.

\subsubsection*{Essentially disjoint unions and essential parts.}
Inquisitive bisimulation between world- or 
state-pointed inquisitive models is robust under the augmentation 
of the set of worlds by disconnected sets of new worlds. This
phenomenon is well known from ordinary bisimulation between Kripke structures.
But whereas the disjoint union of two Kripke models is again a Kripke model,
the disjoint union of two relational inquisitive models would fail to satisfy
extensionality (and thus fail to be a relational inquisitive model)
unless we take care of identifying the respective empty states. 
In the following, if $\str{M}$ and $\str{M'}$ are relational models,
we denote by $\str{M}\oplus\str{M'}$ their \emph{essentially disjoint union},
i.e., the model obtained from the disjoint union by identifying the
empty information states of $\str{M}$ and $\str{M}'$.
Indeed, the the empty information state plays a
special, albeit somewhat trivial r\^ole for various purposes.
To isolate the structurally distinctive part of a 
relational inquisitive model $\str{M} = ( W, S, E, \in, (P_i))$ 
we may consider its \emph{essential part} as obtained by removal of 
the empty information state:
\[
\str{M}^\circ := \str{M}\restr  ( W \cup S^\circ)
= ( W, S^\circ, E^\circ , {\in}^\circ, (P_i)) 
\]
where $S^\circ := S \setminus \{ \emptyset \}$ and 
$E^\circ,{\in}^\circ \subset W \times S^\circ$ are corresponding restrictions.
While $\str{M}^\circ$ is not itself a relational inquisitive model, it
uniquely determines the original model. Writing 
$\str{M}^\circ \ast \{ \emptyset \}$ for the unique extension by re-insertion 
of $\emptyset$, which reproduces $\str{M}$, we have the one-to-one
correspondence 
\[
(\dagger) \qquad
\str{M}^\circ= \str{M}\restr  ( W \cup S^\circ)
\quad \leftrightsquigarrow \quad
\str{M} = \str{M}^\circ \ast \{ \emptyset \} \quad\;\;\qquad\nt
\]

Moreover, essentially disjoint unions of models (or of subsets of their 
domains) are disjoint unions at the level of the essential parts.

\subsubsection*{Locality and truncation of models.}
Towards the assessment of the expressive power of 
$\FO$ over relevant classes of relational inquisitive models, which are
not elementary, we cannot rely on classical compactness
arguments. Instead we invoke \emph{locality arguments} based on the 
local nature of first-order logic over relational structures,
in terms of \emph{Gaifman distance}. 
In the setting of inquisitive relational models, Gaifman distance is graph distance in the
undirected bi-partite graph on the sets $W$ of worlds and $S$ of
states with edges between any pair linked by $E$ or $\ee$;
the \emph{$\ell$-neighbourhood} $N^\ell(w)$ of a world $w$ consists 
of all worlds or states at distance up to $\ell$ from $w$ in this
sense, and $N^\ell(s)$ is similarly defined.
But the presence of the empty information state $\emptyset \in S$
might seem to spoil any locality-based arguments because it 
trivialises the distance measure in $\str{M}$.%
\footnote{In fact, the Gaifman diameter of any 
relational inquisitive model is easily seen to be bounded by $4$, as 
$\emptyset \in S$ is $E$-related to every world, so that also every information 
state has distance at most~$2$ from $\emptyset$.} 
Passage to the essential part $\str{M}^\circ$, however, overcomes this
obstacle.  The empty state plays a trivial r\^ole not just w.r.t.\ bisimulation,
where it only occurs as a dead end, but also w.r.t.\ $\FO$
expressiveness: the relational model $\str{M}$ 
is uniformly quantifier-free $\FO$-interpretable in
its essential part $\str{M}^\circ$.
It follows that $\equiv_q$ between the essential parts of 
(pointed) relational inquisitive models implies $\equiv_q$ between 
the actual models.%
\footnote{Uniform syntactic rewriting provides for every $\phi(x)
  \in \FO$ a translation
$\phi^\circ(x) \in \FO$ of the same quantifier
  rank such that $\str{M},w \models \phi$ iff $\str{M}^\circ, w \models \phi^\circ$.} 
So the correspondence in~$(\dagger)$ is
compatible with individual levels of $\FO$-equivalence. 
Meaningful locality arguments can therefore be based 
on $\ell$-neighbourhoods w.r.t.\ essential parts, where Gaifman
distance is not trivialised by $\emptyset$.
If $\str{M}$ is a relational inquisitive model, 
$w$ a world in $\str{M}$, and $\ell$ an even number, we define 
\emph{truncation} of $\str{M}$ to depth~$\ell$ as 
\[
\str{M}_w^\ell:=(\str{M}^\circ\restr N^\ell(w))*\{\emptyset\}, 
\]
where $N^\ell(w)$ consists of those worlds or states at Gaifman distance at most $\ell$ 
from $w$ in $\str{M}^\circ$ (!). It is easy to see that $\str{M}_w^\ell$ is a relational 
inquisitive model. Similarly, 
for a state $s \not=\emptyset$ and even $\ell$,
we define $\str{M}_s^\ell:=(\str{M}^\circ\restr N^{\ell+1}(s))*\{\emptyset\}$, 
which again is a relational inquisitive model.

\Section{Characterisation theorem for \inqml}
\label{inqmlcharsec}  
\label{inqbmsec}

Our aim is to show the following main characterisation theorem. 

\medskip
\textsc{Theorem~\ref{main1}}. 
\emph{Let $\CC$ be any of the following classes of relational models: 
the class of all models; of finite models; of locally full models; of 
 finite locally full models. Then
 $\inqbm \equiv_\CC^{\w} \FO/{\sim}$  and  $\,\inqbm \equiv_\CC^{\s} \FO/{\sim}.$}

\medskip
By Observation~\ref{compactobs} 
it suffices to establish the following.

\bP 
\label{compprop}
Let $\CC$ be any of the above classes. Any first-order formula which is $\sim$-invariant over $\CC$ 
is $\sim^n$-invariant over $\CC$ for some finite
level $n \in \N$. 
\eP

\subsection{Partial unfolding and stratification}
\label{stratsec}

To establish the compactness principle for $\sim$-invariance 
expressed in Proposition~\ref{compprop} for the relevant
classes of relational models we make use of a 
process of \emph{stratification}. This is similar to tree-like unfoldings 
in standard modal logic. 

\bD
\label{treelikestratdef}
A relational inquisitive model $\str{M}$ 
is \emph{stratified} if its two domains $W$ and $S$ consist of 
essentially disjoint strata, 
i.e., of two families $(W_i)_{i\in\N}$ and $(S_i)_{i\in\N}$ such that:
\bre
\item
$W ={\bigcup}\, W_{i}$ and
$S = {\bigcup} S_i$;
\item for each $i,j\in\mathbb{N}: W_i\cap W_j=\emptyset$ and $S_i\cap S_j=\{\emptyset\}$;
\item $E[w] \subset S_{i}$
for all $w\in W_{i}$ and $S_i \subset \wp(W_{i+1})$.
\ere
For an even $\ell\neq 0$ and a world $w$, we say that $\str{M}$ 
is \emph{stratified to depth $\ell$ from $w$} if
its truncation $\str{M}_w^\ell$ (see the definition at the end of Section \ref{sec:route})
is stratified
with $W_0 = \{w\}$. 
$\str{M}$ is \emph{stratified to depth $\ell$} from a
non-empty state $s \in S$, if its truncation 
$\str{M}_s^{\ell}$ is
stratified with $W_0 = \emptyset$, $S_0=\wp(s)$. 
\eD

We note that no non-trivial stratified model can be full.

\bP
\label{stratunfoldprop}
Any world-pointed relational inquisitive
model $\str{M},w$ is bisimilar to a stratified one. 
For even $\ell\neq 0$, any finite $\str{M},w$ is bisimilar 
to a finite model that is stratified to depth $\ell$ from $w$.
An analogous result holds for state-pointed relational inquisitive models $\str{M},s$. 
If $\str{M}$ is locally full, the ($\ell$-)stratified target model 
can be chosen to be locally full, too. 
\eP

\prf[Proof Sketch]
The underlying process of partial unfolding is similar 
to the well-known tree unfolding of Kripke structures, but leaves
quite some flexibility as to the choice 
of the second sort. 
Starting from a model $\str{M}$, we define a stratified model $\str{M}'$,
whose essential part $\str{M}'\nt^{\circ}$ consists 
of $\N$-tagged copies 
of worlds and non-empty
information states from $\str{M}$, so that 
$W' \subset W \times \N$ and 
$S'\nt^\circ \subset S^\circ \times \N$.
In the world-pointed case, let $w' := (w,0)$. We take $W'_0 := \{ (w,0) \}$. 
For any $n\in \mathbb{N}$, we choose a downward closed set $S_n\supseteq \bigcup_{(u,n)\in W_n'}E[u] $
and put 
\[
\barr{r@{\;:=\;}l}
S'_n\!\nt^\circ  &  S_n^\circ \times \{n\}, 
\\
W'_{n+1} & \bigcup_{s\in S_n^\circ} s \times \{n+1\}, 
\\
\hnt
\mbox{ and define }
E'\nt^\circ & \{ ((u,n),(s,n))
\colon (u,s) \in E \},
\\
\ee'\nt^\circ & \{((u,n+1),(s,n))\colon u\in s\},
\\
P'_i & \{(u,n)\colon u\in P_i\}. 
\earr
\]
This uniquely determines $\str{M}'=\str{M}'\nt^\circ*\,\{\emptyset\}$.
It is easy to check that $\str{M}',w' \sim \str{M},w$.
In order to maintain finiteness, the unfolding process can be
truncated at any stage $n$ if we replace the above $W'_{n+1}$ by 
$W$ and correspondingly put $S^\circ$ instead of $S'_{n+1}\!\nt^\circ$ and augment 
$E'\nt^\circ$ by all of~$E^\circ$. The resulting $\str{M}',w'$ is fully
bisimilar to $\str{M},w$, is finite if $\str{M}$ is, and is stratified
to depth $2n$. With the straightforward maximal choice for the 
$S'\nt_n^\circ$, viz.\
$S'\nt_n^\circ := S^\circ \times \{n\}$, the (full or truncated) unfolding process
preserves local fullness, too. 
In the state-pointed case we start out by setting $W_0' :=\emptyset$, 
$S_0' := \wp(s)\times \{ 0 \}$ and proceed inductively as above. 
\eprf

\bO
\label{stratcutoffobs} 
Let $\str{M},w$ and $\str{M}',w'$ be world-pointed relational models that are stratified 
to depth $\ell$ for some even $\ell$ 
from their respective worlds. Let $\str{M}_w^\ell,w$ and $\str{M}'\nt_{w'}^\ell,w'$ be their $\ell$-truncations.
Then, for 
$n \geq \ell/2$:
\[
\str{M}_{w}^\ell,w 
\;\simn\; \str{M}'\nt_{w'}^\ell,w'
\;\; \Longrightarrow \;\;
\str{M}_w^\ell,w\,\sim\, \str{M}'\nt_{w'}^\ell,w'.
\]
Analogously for state-pointed models that are stratified to depth
$\ell$ from their distinguished states. 
\eO

This is because, due to stratification and cut-off, the $n$-round game 
exhausts all possibilities in the unbounded game.
After $m$~rounds of the bisimulation game, which starts from the 
pairing of the roots~$w$ and $w'$ of the stratified models, 
the position is a pairing of worlds from strata $W_m$ and $W_m'$.
So player~\PlayerII\ wins the unbounded game if she does not lose within 
$n$~rounds.

\prf[Proof of Theorem~\ref{main1}.]
We first present the upgrading
argument for the case of world-pointed models,
which is closer to the classical intuition. 
The version for state-pointed models, which is formally the stronger, 
will be discussed below.
Let $\CC$ be any one of the classes in the theorem 
and let $\phi(x) \in \FO_q$ be $\sim$-invariant as a world property 
over $\CC$. We want to show that $\phi$ is $\simn$-invariant 
over $\CC$ for $n = 2^q$, where $q$ is the quantifier rank of~$\phi$. 
The upgrading argument is sketched in Figure~\ref{iqbmfigure}.
Towards its ingredients, consider 
a world-pointed relational model $\str{M},w$ 
in $\mathcal{C}$. Since $\phi$ is $\sim$-invariant, 
we can, by Proposition~\ref{stratunfoldprop}, 
assume w.l.o.g.\ that $\str{M},w$ is stratified to depth $\ell = n$ from $w$. 
Let $\str{M}_w^\ell$ be its $\ell$-truncation, which is then fully stratified. 
We define two world-pointed models $\str{M}_0,w$ and $\str{M}_1,w$ as
follows. 
Each of the  $\str{M}_i$ consists of an essentially disjoint union
of the following constituents: 
both models contain $q$ distinct isomorphic copies
of $\str{M}$ as well as of $\str{M}_w^\ell$. 
In addition, $\str{M}_0$ contains a copy 
of 
$\str{M}_w^\ell$
with the distinguished world $w$,
 while $\str{M}_1$ contains a copy of $\str{M}$ with the 
 distinguished world~$w$.%
\[
\barr{@{}r@{\;:=\;\;}r@{\;\;\oplus\;\;}c@{\;\;\oplus\;\;}l@{}}
\str{M}_0,w &
q \otimes \str{M} &
\str{M}_w^\ell,w &
q \otimes \str{M}_w^\ell
\\
\str{M}_1,w &
q \otimes \str{M} &
\str{M},w &
q \otimes \str{M}_w^\ell
\earr
\]

Using a locality-based Ehrenfeucht-Fra\"iss\'e game argument for $\FO$ 
we can show:
\[
(\ast) \quad \str{M}_0,w  \equiv_q \str{M}_1,w.
\]

As argued in connection with the correspondence~$(\dagger)$ at the end
of Section~\ref{bisiminvsec}, due to quantifier-free interpretability
of $\str{M}_i,w$ in $\str{M}_i^\circ,w$, 
$(\ast)$ is equivalent to 
\[
(\ast\ast) \quad \str{M}_0^{\circ},w  \equiv_q
\str{M}_1^{\circ},w 
\]

The diagram in Figure~\ref{disjointconesfigure}
suggests the arrangement, with open cones for 
copies of $\str{M}^\circ$ and truncated cones for 
$\str{M}^\circ \restr N^\ell(w)$ (the essential part of $\str{M}_w^\ell$),
and with filled circles for the distinguished worlds.

\begin{figure}  
\[
\makebox(100,110){$
\underbrace{
\makebox(60,100)[b]{$
\begin{xy}
(0,0) *{\circ} \ar@{-},(-8,20),
+(0,0) \ar@{-},(8,20),
\end{xy}
\quad
\begin{xy}
(0,0) *{\circ} \ar@{-},(-8,20),
+(0,0) \ar@{-},(8,20),
\end{xy}
$}}_{\mbox{\tiny $q$ copies}}
\quad\quad
\begin{xy}
(0,-3)*{w},
(0,0) *{\bullet},(-4,10) **@{-},
(0,0) ,(4,10) **@{-},
(-4,10);(4,10) **@{-}
\end{xy}
\quad\;
\underbrace{
\makebox(40,100)[b]{$
\begin{xy}
(0,0) *{\circ};(-4,10) **@{-},
(0,0) *{\circ};(4,10) **@{-},
(-4,10);(4,10) **@{-}
\end{xy}
\quad
\begin{xy}
(0,0) *{\circ};(-4,10) **@{-},
(0,0) *{\circ};(4,10) **@{-},
(-4,10);(4,10) **@{-}
\end{xy}
$}}_{\mbox{\tiny $q$ copies}}
\qquad
\equiv_q
\qquad
\underbrace{
\makebox(60,100)[b]{$
\begin{xy}
(0,0) *{\circ} \ar@{-},(-8,20),
+(0,0) \ar@{-},(8,20),
\end{xy}
\quad
\begin{xy}
(0,0) *{\circ} \ar@{-},(-8,20),
+(0,0) \ar@{-},(8,20),
\end{xy}
$}}_{\mbox{\tiny $q$ copies}}
\quad\quad\quad
\begin{xy}
(0,-3)*{w},
(0,0) *{\bullet} \ar@{-},(-8,20),
+(0,0) \ar@{-},(8,20)
\end{xy}
\quad
\underbrace{
\makebox(40,100)[b]{$
\begin{xy}
(0,0) *{\circ};(-4,10) **@{-},
(0,0) *{\circ};(4,10) **@{-},
(-4,10);(4,10) **@{-}
\end{xy}
\quad
\begin{xy}
(0,0) *{\circ};(-4,10) **@{-},
(0,0) *{\circ};(4,10) **@{-},
(-4,10);(4,10) **@{-}
\end{xy}
$}}_{\mbox{\tiny $q$ copies}}
$
}
\]
\caption{The structures $\str{M}_0^\circ,w$ and $\str{M}_1^\circ,w$ in the game argument.}
\label{disjointconesfigure}
\end{figure}

We argue that the second player has
a winning strategy in the classical $q$-round Ehrenfeucht--Fra\"\i ss\'e game 
over the two structures in $(\ast\ast)$ starting in the position with a 
single pebble on the distinguished world $w$ on either side. 
Indeed, player \PlayerII\ can force a win by maintaining
the following invariant w.r.t.\ the game positions 
$(\ubar;\ubar')$ for 
$\ubar = (u_0,u_1,\ldots, u_m)$ with $u_0 = w$ in $\str{M}_0^\circ$ and 
$\ubar' = (u'_0,u_1',\ldots, u_m')$ with $u'_0 = w$ in $\str{M}_1^\circ$ 
after round $m$, for $m=0,\ldots, q$, for $\ell_m := 2^{q-m}$:

\begin{quotation}\noindent
$\ubar$ and $\ubar'$ are partitioned into clusters of matching
sub-tuples such that the distance between separate clusters is greater
than $\ell_m$ and corresponding clusters are in isomorphic
configurations of isomorphic component structures of $\str{M}_0^\circ$ and
$\str{M}_1^\circ$ or in isomorphic configurations in
$\str{M}_0^\circ\restr N^\ell(w)$ and $\str{M}_1^\circ\restr N^\ell(w)$. 
\end{quotation}

This condition is satisfied at the start of the game, for $m=0$
($\ell_0 = 2^q = n$).
The second player can maintain this condition 
through a round, say in the step from $m$ to $m+1$, as follows.
Suppose the first player puts a pebble in position $u= u_{m+1}$ in
$\str{M}_0^\circ$ or $u' = u'_{m+1}$ in $\str{M}_1^\circ$ at distance up to $\ell_{m+1}$ 
of one of the level $m$ clusters (it cannot fall within  distance
$\ell_{m+1}$ of two distinct clusters, since the distance between two
distinct clusters from the previous level is greater than $\ell_m = 2
\ell_{m+1}$); then this new position joins a sub-cluster of that cluster 
and its match is found in an isomorphic position relative to the matching cluster.
If the first player puts the new pebble in a position $u= u_{m+1}$ in
$\str{M}_0^\circ$ or $u' = u'_{m+1}$ in $\str{M}_1^\circ$ at distance greater 
than $\ell_{m+1}$ of each one of the level $m$ clusters, 
this position will form a new cluster and 
can be matched with an isomorphic position in one of the
as yet unused component structures on the opposite~side.

This argument restricts naturally to the
scenarios of (finite or general) locally full relational
inquisitive structures, because stratification (to some depth)
according to Proposition~\ref{stratunfoldprop} 
preserves local fullness,
and so does restriction to some even depth and the formation of essentially
disjoint sums. 

Given any two pointed models $\str{M},w\,\simn\,\str{M}',w'$ in any of
the relevant classes $\CC$, we see that a first-order formula
$\phi$ of quantifier rank $q$ that is preserved under $\sim$,
is preserved by chasing the diagram in Figure~\ref{iqbmfigure}
along the path through the auxiliary models,  
which are all in $\mathcal{C}$.
The expressive completeness claim for Theorem~\ref{main1}, i.e.\ expressibility
of $\phi$ in $\inqml$ over $\CC$, now follows 
from Corollary~\ref{EFcorrworldpointedrel}: indeed, $\phi$ is logically
equivalent over $\CC$ to the disjunction over the characteristic
formulae  $\chi^n_{\M,w}$ for all $\str{M},w \in \CC$ that satisfy $\phi$.
\eprf

\begin{figure}
\[
\xymatrix{
\str{M},w \ar@{-}[d]|{\rule{0pt}{1ex}\sim} \ar @{-}[rr]|{\;\simn\,}
&& \str{M}',w' \ar @{-}[d]|{\rule{0pt}{1ex}\sim}
\\
\str{M}_1,w \ar@{-}[d]|{\rule{0pt}{1.2ex}\;\equiv_q}
\ar@{.}[rr]|{\;\simn\,}
&& \str{M}'_1,w'
\ar@{-}[d]|{\rule{0pt}{1.2ex}\;\equiv_q}
\\
\str{M}_0,w \ar@{-}[d]|{\rule{0pt}{1ex}\sim}
\ar@{.}[rr]|{\;\simn\,}
&& \str{M}'_0,w'
\ar@{-}[d]|{\rule{0pt}{1ex}\sim}
\\
\str{M}_w^\ell,w \ar@{-}[rr]|{\;\sim\,}
&& 
\str{M}'_{w'}\!\!\!\nt^\ell,w'
}
\]
\caption{Upgrading pattern for Theorem~\ref{main1}/Proposition~\ref{compprop}.}
\label{iqbmfigure}
\end{figure}

\subsubsection*{The case of state properties.}
To show Proposition~\ref{compprop} for state properties, 
we can similarly upgrade the situation $\str{M},s \,\simn\,
\str{M}',s'$ 
for non-empty $s,s'$ in companion structures through passage to truncations of 
fully bisimilar models that are stratified to depth $\ell$ from 
their distinguished states.%
\footnote{Note that the upgrading argument trivialises for $\phi(\s) \in \FO$ 
in the case of state-pointed models $\str{M},s$ with $s = \emptyset$: 
For any $n$, $\str{M},\emptyset \,\,\simn\,\, \str{M}',s'$ iff
$\str{M},\emptyset \sim \str{M}',s'$ iff $s'=\emptyset$.  
If $\phi(\s)$ is $\sim$-invariant over $\CC$, it must be satisfied by 
$s = \emptyset$ across all of $\CC$ or nowhere; 
and being inquisitive over $\CC$, the former must in fact be true.} 
Assuming w.l.o.g.\ that $\str{M},s$ is itself~stratified to depth 
$\ell = 2^q$ with stratified restriction $\str{M}_s^\ell$  
we define as before the following essentially disjoint unions 
\[
\barr{@{}r@{\;:=\;\;}r@{\;\;\oplus\;\;}c@{\;\;\oplus\;\;}l@{}}
\str{M}_0,s &
q \otimes \str{M} &
\str{M}_s^\ell,s &
q \otimes \str{M}_s^\ell
\\
\str{M}_1,s &
q \otimes \str{M} &
\str{M},s &
q \otimes \str{M}_s^\ell
\earr
\]
and we find that $\str{M}_0,s  \equiv_q \str{M}_1,s$. 
We do the same for $\str{M}',s'$. 
The rest of the argument for Proposition~\ref{compprop}  
is completed with the straightforward analogue of Figure~\ref{iqbmfigure}
for the relevant state-pointed models.

\section{Conclusion}

We have seen the beginnings of a model theory for
inquisitive modal logic.
Our contribution started in Section~\ref{sec:bisimulations}, 
where we described a natural notion of bisimulation for 
inquisitive modal structures. 
From a game-theoretic perspective, bisimilarity and its 
approximations can be characterised in terms of a game which 
interleaves two kinds of moves: world-to-state moves 
(from $w$ to some $s\in\Sigma(w)$) and state-to-world moves 
(from $s$ to some $w\in s$).

In Section~\ref{sec:ef} we saw that bisimilarity relates to modal 
equivalence in the usual way: two pointed models 
over a finite vocabulary 
are distinguishable in the $n$-round bisimulation game 
iff they are distinguished by a formula of modal depth~$n$.

In Section~\ref{sec:neighbourhood} we compared inquisitive modal logic 
to neighbourhood semantics for modal logic, showing that, although these 
two logics are interpreted over similar structures, they are very different in 
terms of their expressive power, and are invariant under different notions 
of bisimulation equivalence. 

In Section~\ref{relinqmodsec} we discussed how inquisitive modal models 
can be encoded as two-sorted relational structures on which we can 
naturally interpret first-order formulae of a suitable relational signature. 
This enabled us to define a standard translation from \inqbm\ to 
first-order logic, and to view \inqbm\ as a syntactic fragment of 
first-order logic with respect to those relational structures.

We then asked over what classes of structures this syntactic fragment coincides, up to logical equivalence, with the 
fragment determined by the semantic property of bisimulation invariance.
Using an inquisitive analogue of partial tree unfoldings in
Section~\ref{inqbmsec}, we established a positive answer to 
this question for several natural classes, including the class of all
relational inquisitive models, and the class of all finite models.

The results obtained in this paper provide us with a better 
understanding of inquisitive modal logic in at least two ways.
From a more concrete perspective, we have given a characterisation 
of the expressive power of \inqbm\ which is very helpful in order 
to tell what properties of pointed models can and cannot be expressed 
in the language: for instance, it is easy to see that properties like 
$\mathcal{P}(w):=``W \in\Sigma(w)\text{''}$ or 
$\mathcal{P}(w):=``\{w\}\in\Sigma(w)\text{''}$ are not bisimulation
invariant, and thus not expressible in \inqml.
From a more abstract perspective, we have looked at a natural notion of 
behavioural equivalence for inquisitive modal structures, whose main 
constituent is a map $\Sigma:W\to\wp\wp(W)$, rather than
$\sigma:W\to\wp(W)$ 
as in Kripke structures. We saw that, in terms of expressive power, 
\inqml\ is a natural choice for a language designed to talk about
properties which are invariant under this notion of equivalence: 
over various classes of structures, \inqbm\ expresses all and only 
the first-order properties that are invariant in this sense.

In a separate paper, we will tackle the case of inquisitive
epistemic models---the inquisitive version of multi-modal $S5$
models. Conceptually, this class of models is interesting in light of the 
natural interpretation of inquisitive modalities in the epistemic setting, 
as described in Section~\ref{sec:inquisitive modal logic}. Technically, it presents
interesting challenges as the stratifications 
used in Section~\ref{inqmlcharsec} are incompatible with the 
$S5$ frame conditions. Nevertheless, it can be shown that the counterpart 
of our characterisation result still holds in this setting---again, both in general and in restriction to finite models.
A preliminary account is given in~\cite{CiardelliOttoarXiv18}.

\bibliographystyle{asl}

\end{document}